\documentclass[12pt]{article}
\usepackage{amssymb,amsfonts,amsmath,amsthm,indentfirst}
\renewcommand{\theequation}{\arabic{section}.\arabic{equation}}

\newtheorem{theorem}{Theorem}[section]
\newtheorem{lemma}[theorem]{Lemma}
\newtheorem{remark}{Remark}[section]

\title{
\Large Singular limit and long-time dynamics of Bresse systems
}
\author{
\small To Fu Ma \thanks{Email: matofu@icmc.usp.br. (Corresponding Author)} \\
\small Rodrigo Nunes Monteiro \thanks{Email: rodrigonunesmonteiro@gmail.com.} \\
\small Institute of Mathematical and Computer Sciences, University of S\~ao Paulo \\
\small 13566-590 S\~ao Carlos, SP, Brazil
}

\date{}

\begin{document}

\maketitle

\begin{abstract}
The Bresse system is a valid model for arched beams which reduces to the classical Timoshenko system when the arch 
curvature $\ell=0$. Our first result shows the Timoshenko 
system as a singular limit of the Bresse system as $\ell \to 0$. 
The remaining results are concerned with the long-time dynamics of Bresse systems. 
In a general framework, allowing nonlinear damping and forcing terms, 
we prove the existence of a smooth global attractor with finite fractal dimension   
and exponential attractors as well. 
We also compare the Bresse system with the Timoshenko system, 
in the sense of the upper-semicontinuity of their attractors as $\ell \to 0$.
\end{abstract}



\setcounter{equation}{0}
\section{Introduction}
In this paper we study the long-time dynamics of solutions of a semilinear Bresse system for vibrations of curved beams.
The linear system is given by three motion
equations,
\begin{align}
 \rho_{1}\varphi_{tt} - k(\varphi_{x}+\psi+\ell w)_{x} - k_{0}\ell (w_{x}-\ell\varphi) = 0 & \;\;  \mbox{in} \:\: (0,L) \times \mathbb{R}^{+} ,  \label{BO1}\\
 \rho_{2}\psi_{tt} - b\psi_{xx}+k(\varphi_{x} + \psi+\ell w) = 0  & \;\;  \mbox{in} \: \:(0,L) \times \mathbb{R}^{+},  \label{BO2} \\
 \rho_{1}w_{tt} - k_{0}(w_{x}-\ell\varphi)_{x} + k\ell(\varphi_{x}+\psi+\ell w) =0 & \;\; \mbox{in} \: \:(0,L) \times \mathbb{R}^{+}, \label{BO3}
\end{align}
where $\varphi, \psi, w$ represent, respectively, vertical displacement, shear angle, and longitudinal displacement. The coefficients
$\rho_1, \rho_2, b , k, k_0$ are positive constants related to the material and the parameter $\ell$ stands for the curvature of the beam.
In the context of a circular arch of radius $R$ one has $\ell=R^{-1}$. A description of the model can be found in \cite[Chap. 3]{lagnese-book}.
The original derivation of Bresse system was presented in \cite{Bresse}.

\bigskip

We observe that when the curvature $\ell \to 0$ the Bresse system uncouples into the well-known Timoshenko
system,
\begin{align}
\rho_{1}\varphi_{tt} - k(\varphi_{x}+\psi)_{x} =0 &  \;\; \mbox{in} \:\:(0,L) \times \mathbb{R}^{+},  \label{TO1} \\
\rho_{2}\psi_{tt} - b\psi_{xx}+k(\varphi_{x} + \psi) =  0  & \;\; \mbox{in} \: \:(0,L) \times \mathbb{R}^{+}, \label{TO2}
\end{align}
and an independent wave equation $\rho_{1} w_{tt}-k_0 w_{xx}=0$.
Therefore sometimes the Timoshenko system is called Bresse-Timoshenko system. The derivation of the Timoshenko
system is presented in \cite{Timoshenko}.

There are several works dedicated to the mathematical analysis of the Bresse system.
They are mainly concerned with decay rates of solutions of the linear system. This is done
by adding suitable damping effects that can be of thermal, viscous or  viscoelastic nature.
By analogy to the Timoshenko system, a remarkable stability criteria for the Bresse system is the equal wave speeds assumption
\begin{equation} \label{equal-speed}
\frac{\rho_1}{k} = \frac{\rho_2}{b} \quad \mbox{and} \quad k=k_0.
\end{equation}
Roughly speaking, if damping terms are added only in one or two of the
equations in the Bresse system then exponential stability is only possible if the equal wave speeds assumption 
(\ref{equal-speed}) holds. Of course, if weak damping are added in all the equations of the system then
assumption (\ref{equal-speed}) is not needed for exponential stability. See for instance
\cite{alabau,alves,delloro,fatori-monteiro,fatori-munoz,liu-rao,santos-soufyane,soriano-jmaa,soufyane-said,wehbe}.

\medskip

On the other hand, it is worthy mentioning that all above cited works on Bresse systems were concerned with linear problems.
With respect to nonlinear problems, the stability of the Bresse system was studied in \cite{charles-soriano}, 
with three locally defined nonlinear damping without assuming the equal speeds assumption. 
An analogous result for Timoshenko systems was presented in \cite{cavalcanti-timoshenko}.

It turns out that long-time dynamics characterized by global attractors 
was not discussed for Bresse systems. 
Even for the Timoshenko system there are only a few works in this direction. 
Here we consider the nonlinear model,
\begin{align}
\rho_{1}\varphi_{tt} - k(\varphi_{x}+\psi+\ell w)_{x} - k_{0}\ell (w_{x}-\ell\varphi)+g_{1}(\varphi_t) + f_1(\varphi,\psi,w) =0, & \label{P01} \\
\rho_{2}\psi_{tt} - b\psi_{xx}+k(\varphi_{x} + \psi+\ell w) + g_{2}(\psi_t)+ f_2(\varphi,\psi,w) = 0, &  \\
\rho_{1}w_{tt} - k_{0}(w_{x}-\ell\varphi)_{x} + k\ell(\varphi_{x}+\psi+\ell w) + g_{3}(w_t) + f_3(\varphi,\psi,w) =0, &  \label{P03}
\end{align}
defined in $(0,L) \times \mathbb{R}^{+}$, where $g_{1}(\varphi_t) ,g_{2}(\psi_t), g_{3}(w_t)$ are nonlinear damping terms
and $f_i(\varphi,\psi,w)$, $i=1,2,3$, are nonlinear forces.
To this system we add Dirichlet boundary
condition
\begin{equation} \label{bci}
\varphi(0,t) = \varphi(L,t) = \psi(0,t) = \psi(L,t) = w(0,t) = w(L,t) = 0, \quad t \in \mathbb{R}^{+},
\end{equation}
and initial condition
\begin{equation} \label{3ci}
\varphi(0) = \varphi_{0},\:\: \varphi_t(0) = \varphi_{1},\:\: \psi(0) = \psi_{0},\:\: \psi_t(0) = \psi_{1}, \:\: w(0)= w_0, \:\: w_t(0) = w_1.
\end{equation}
Since our problem has damping terms in all of the equations (\ref{P01})-(\ref{P03}) we shall not assume the equal wave speeds assumption.

The main features of the paper are summarized as follows.

\bigskip

$(i)$  We present a rigorous proof showing 
that solutions of the Bresse system converge to that of the Timoshenko system 
as $\ell \to 0$. This is called singular limit because at $\ell = 0$ the Bresse system uncouples. 
See Theorem \ref{theo-singular}.

\bigskip 

$(ii)$ By considering a nonlinear damping, without growth restriction near  zero,  
we establish the existence of a global attractor. We also show that the system is gradient and therefore the attractors
are characterized as unstable manifold of the set of stationary solutions. 
See Theorem \ref{theo-attractor1}. 

\bigskip

$(iii)$ By assuming further that damping terms are Lipschitz, we derive a stability inequality and
prove that the global attractor is regular and has finite fractal dimension. The same hypotheses 
allow us to show the existence of generalized exponential attractors as well. See Theorem \ref{theo-attractor2}.

\bigskip

$(iv)$ We also compare the attractors of the Bresse system with those of the Timoshenko system. 
More precisely, we show the upper semicontinuity of attractors of (\ref{P01})-(\ref{3ci}) as $\ell \to 0$. 
In this case we shall assume that $f_1,f_2$ are not depending on $w$. This is reasonable since in the limit
$\ell = 0$ we obtain the Timoshenko system, where longitudinal displacement $w$ is neglected.
That is, after limit, we get
\begin{align} \label{TS101}
 \rho_{1}\varphi_{tt} - k(\varphi_{x} + \psi)_{x} + g_{1}(\varphi_t) + f_1(\varphi,\psi)  = 0 & \;\; \mbox{in} \:\: (0,L) \times \mathbb{R}^{+} , \\
 \rho_{2}\psi_{tt} - b\psi_{xx}+k(\varphi_{x} + \psi ) + g_{2}(\psi_t)+ f_2(\varphi,\psi)  = 0 & \;\;  \mbox{in} \: \:(0,L) \times \mathbb{R}^{+},
\end{align}
subjected to Dirichlet boundary
condition
\begin{equation} \label{TS102}
\varphi(0,t) = \varphi(L,t) = \psi(0,t) = \psi(L,t) =0 , \quad  t \in \mathbb{R}^{+},
\end{equation}
and initial condition
\begin{equation} \label{TS103}
\varphi(0) = \varphi_{0},\:\: \varphi_t(0) = \varphi_{1},\:\: \psi(0)= \psi_{0},\:\: \psi_t(0) = \psi_{1}.
\end{equation}
See Theorem \ref{theo-upper}.

\bigskip

\noindent $(v)$ Above, one of the difficulty is to obtain uniform estimates for global attractors 
of the Bresse system with respect to the curvature $\ell$. This is done by showing the existence of an absorbing set that is uniformly bounded with respect to $\ell$ (small). 
See Lemma \ref{lemma-AbsobSet}.  

\setcounter{equation}{0}
\section{Preliminaries and well-posedness}
We begin with some notations on the standard $L^p(0,L)$ and Sobolev $H^m(0,L)$ spaces. The $L^p$ norm is denoted by
$$
\| u \|_p \;\;  \mbox{if} \;\; p \neq 2 \quad \hbox{and} \quad \| u \| \;\; \mbox{if} \;\; p=2 .
$$
For the Sobolev space $H^1_0(0,L)$ we have  
$$
\Vert u \Vert \leq \frac{L}{\pi} \Vert u_x \Vert
\quad \hbox{and} \quad
\Vert u \Vert_{H^1_0} = \Vert u_x \Vert. 
$$
The phase space we consider is that of weak solutions
\begin{equation}
\mathcal{H} = H^1_0 (0,L)^3 \times L^2 (0,L)^3,
\end{equation} 
equipped with the standard norm 
$$
\| y \|_{\mathcal{H}}^2 =  \|\varphi_x \|^2 + \|\psi_x\|^2 + \|w_x \|^2 + \|\tilde{\varphi}\|^2 + \|\tilde{\psi}\|^2+ \|\tilde{w}\|^2 ,
$$ 
where $y=(\varphi, \psi, w, \tilde{\varphi}, \tilde{\psi}, \tilde{w})$. For convenience, sometimes we use the ($\ell$-dependent) norm,
\begin{equation} \label{l-norm}
\| y \|_{\mathcal{H}_{\ell}}^2 =
\rho_1\|\tilde{\varphi}\|^2 +\rho_2\|\tilde{\psi}\|^2+\rho_1\|\tilde{w}\|^2 + b\|\psi_x\|^2 + k \|\varphi_x + \psi +\ell w\|^2+ k_0 \|w_x - \ell \varphi \|^2,
\end{equation}
which is equivalent to the standard norm. 
Indeed, clearly there exists $\gamma_1>0$ such that
$\| y \|_{\mathcal{H}_{\ell}}^2 \leq  \gamma_1 \| y \|_{\mathcal{H}}^2$.
Then from the open mapping theorem, there exists $\gamma_2>0$ such
\begin{equation} \label{norma-2}
\| y \|_{\mathcal{H}}^2 \leq  \gamma_2  \| y \|_{\mathcal{H}_{\ell}}^2,
\end{equation}
which shows the equivalence of the norms. In particular there exists
$\gamma_3 >0$ such that
\begin{equation}  \label{norma-ok}
\|\varphi_x \|^2 + \|\psi_x\|^2 + \|w_x \|^2 \leq  \gamma_3  \big( b\|\psi_x\|^2 + k \|\varphi_x + \psi +\ell w\|^2+ k_0 \|w_x - \ell \varphi \|^2 \big).
\end{equation}

\begin{remark} \label{rem-ell} \rm In the study of continuity of attractors as $\ell \to 0$, 
some energy estimates, uniform with respect to $\ell$, are used.  
To this end we need $\gamma_1,\gamma_2,\gamma_3$ independent of $\ell$, for $\ell$ small. 
It is clear that we can choose such $\gamma_1$ 
if $\ell \leq \ell_0$, for some $\ell_0$. 
Here, we show a simple argument to obtain $\gamma_3$ independently $\ell \in [0,\ell_0]$ with $\ell_0 < \pi/2L$. Given $\varphi,\psi,w \in H^1_0(0,L)$, 
\begin{align*}
\| \varphi_x\|^2 & + \| \psi_x\|^2 + \| w_{x }\|^2 \\ 
& \le 
\|\psi_x\|^2 +  2 \| \varphi_x + \psi + \ell w \|^2 
+ 2 \|w_x- \ell\varphi\|^2   
+ 4 \| \psi \|^2 + 4 \ell^2 \| w \|^2 + 2 \ell^2 \|\varphi\|^2 \\
& \le 
\Big(1+ \frac{4L^2}{\pi^2}\Big)\| \psi_x\|^2 + 2 \| \varphi_x + \psi + \ell w \|^2 
+ 2 \|w_x- \ell\varphi\|^2  
\\ 
& \quad \, + \frac{4\ell^2 L^2}{\pi^2} \| w_x\|^2 +
\frac{2 \ell^2 L^2}{\pi^2} \| \varphi_x\|^2.
\end{align*}
Then, for $\ell \in [0, \ell_0]$,  
\begin{align*}
\| \varphi_x\|^2 & + \| \psi_x\|^2 + \| w_{x }\|^2 \\
& \le \Big(1 - \frac{4\ell_0^2 L^2}{\pi^2}\Big)^{-1} 
\Big( 
\Big(1+\frac{4L^2}{\pi^2}\Big) 
\|\psi_x\|^2 +  2 \| \varphi_x + \psi + \ell w \|^2 
+ 2 \|w_x- \ell\varphi\|^2  \Big).
\end{align*}
Hence there exists a constant $\gamma_3>0$ such that  $(\ref{norma-ok})$ holds for all $\ell \in [0,\ell_0]$. 
In this case, we take $\gamma_2 = \max\{ (\min\{ \rho_1,\rho_2\})^{-1}, \gamma_3 \}$. \qed 
\end{remark}

\subsection{Assumptions}

The assumptions we make on $f_1,f_2,f_3$ are those of locally Lipschitz
and gradient type.
Let us assume there exists a $C^2$ function $F: \mathbb{R}^3 \to \mathbb{R}$
such that
\begin{equation} \label{new-F0}
\nabla F = (f_1,f_2,f_3),
\end{equation}
and satisfies the following conditions: There exist $\beta , m_F \ge 0$
such that
\begin{equation} \label{new-F2}
F(u,v,w) \geq  - \beta \big( |u|^{2} + |v|^{2} + |w|^2 \big) - m_F  , \quad \forall \, u,v,w \in \mathbb{R},
\end{equation}
where
\begin{equation} \label{Fc}
0\leq \beta < \frac{\pi^2}{2 \gamma_3 L^2},
\end{equation}
and there exist $p \geq 1$ and $C_f > 0$ such that,
for $i=1,2,3$,
\begin{equation} \label{new-f1}
| \nabla f_i(u,v,w) | \leq
C_f \big( 1 + |u|^{p-1} + |v|^{p-1} + |w|^{p-1} \big), \quad \forall \, u,v,w \in \mathbb{R}.
\end{equation}
In particular this implies that there exists $C_F>0$
such that
\begin{align*} \label{new-F1}
F(u,v,w) \leq C_F \big( 1 + |u|^{p+1} + |v|^{p+1} + |w|^{p+1} \big)  , \quad \forall \, u,v,w \in \mathbb{R}.
\end{align*}
Furthermore, we assume that, for all $u,v,w \in \mathbb{R}$,
\begin{equation} \label{new-fF}
\nabla F(u,v,w)\cdot (u,v,w) - F(u,v,w) \geq
- \beta \big( |u|^{2} + |v|^{2} + |w|^{2} \big) - m_F .
\end{equation}

\begin{remark} \rm A simple example of $F$ satisfying all above assumptions is
$$
F(u,v,w)= |u+v|^4 - |u+v|^2 + \alpha_1 |uv|^2 + \alpha_2 |w|^3, \;\;
\alpha_1 , \alpha_2 \ge 0.
$$
In this case we have
$$
F(u,v,w) \ge \min_{z\in \mathbb{R}}\{z^4-z^2\} = -\frac{1}{4} ,
$$
and
\begin{align*}
&& f_1(u,v,w) = \frac{\partial F}{\partial u} = 4 (u+v)^3 - 2 (u+v)+ 2 \alpha_1 uv^2, \\
&& f_2(u,v,w) = \frac{\partial F}{\partial v} = 4 (u+v)^3 - 2 (u+v)+ 2 \alpha_1 u^2v, \\
&& f_3(u,v,w) = \frac{\partial F}{\partial w} = 3 \alpha_2 |w|w.
\end{align*}
Then conditions $(\ref{new-F2})$-$(\ref{new-f1})$ hold with $m_F= 1/4$ and $p=3$. 
In addition,
$$
\nabla F(u,v,w) \cdot (u,v,w) - F(u,v,w)
\geq 3|u+v|^4 - |u+v|^2
\geq - \frac{1}{16} \ge -m_F,
$$
which shows that $(\ref{new-fF})$ also holds. \qed
\end{remark}

With respect to the damping functions $g_{i} \!\in\! C^1(\mathbb{R})$, $i=1,2,3$, we assume that
\begin{equation}\label{HG1}
g_i \mbox{ is increasing and } g_{i}(0)= 0,
\end{equation}
and there exist constants $m_{i}, M_{i}>0$ such that
\begin{equation} \label{HG2}
m_{i} \, s^2 \leq g_{i}(s)s \leq M_{i} \,s^2, \quad \forall \, |s|> 1.
\end{equation}
To establish the regularity and finite dimension of the attractors
we assume further that (\ref{HG2}) holds for all $s \in \mathbb{R}$, that is
\begin{equation} \label{HG3}
m_{i} \leq g^{\prime}_{i}(s) \leq M_{i}, \quad \forall \, s \in \mathbb{R}.
\end{equation}

\begin{remark} \rm We observe that conditions $(\ref{HG1})$ and $(\ref{HG2})$ imply that for
any given $\delta>0$ there exists $C_{\delta}>0$ such
that
\begin{equation} \label{HG4}
C_{\delta} (g_{i}(u)-g_{i}(v))(u - v) + \delta \geq |u- v|^{2},
\quad \forall \,  u,v \in \mathbb{R}, 
\end{equation}
cf. \cite[Proposition B.1.2]{Yellow}. On the other hand,
assumption $(\ref{HG3})$
implies promptly the usual monotonicity property
\begin{equation} \label{HG5}
(g_{i}(u)-g_{i}(v)) (u - v) \geq  m_i |u- v|^{2},
\quad \forall \, u,v \in \mathbb{R}.  
\end{equation}
\end{remark}

\subsection{Energy identities}
The linear energy of the system, along a solution $(\varphi,\psi,w)$, 
is defined by
%
%
\begin{equation} \label{energy-linear}
E_{\ell}(t) 
 = \frac{1}{2} \Vert (\varphi(t), \psi(t), w(t),\varphi_t(t), \psi_t(t), w_t(t)) \Vert_{\mathcal{H}_{\ell}}^2,
\end{equation}
where $\| \cdot \|_{\mathcal{H}_{\ell}}$ is defined in (\ref{l-norm}). Adding 
forcing terms it becomes
\begin{equation} \label{energy-def}
\mathcal{E}_{\ell}(t) = E_{\ell}(t) + \int_{0}^{L} F(\varphi,\psi,w)\, dx .
\end{equation}
Then, multiplying (\ref{P01})-(\ref{P03}) by $\varphi,\psi,w$ respectively, we obtain 
(formally) by integration over $[0,L]$,  
$$
\frac{d}{dt} \mathcal{E}_{\ell}(t) = - \int_{0}^{L} \big( g_{1}(\varphi_t)\varphi_{t}+g_{2}(\psi_t)\psi_{t}+g_{3}(w_t)w_{t} \big) \, dx, \;\;  t \ge 0,  
$$
which holds for strong solutions. Integrating with respect to $t$ we obtain  
the energy identity  
\begin{equation} \label{energy-ID}
 \mathcal{E}_{\ell}(t) + \int_{s}^{t}\int_{0}^{L}g_{1}(\varphi_t)\varphi_{t}+g_{2}(\psi_t)\psi_{t}+g_{3}(w_t)w_{t} \, dx d\tau =  \mathcal{E}_{\ell}(s) , \;\; 0 \le s < t.
\end{equation}
Moreover, we have the following useful energy inequality.

\begin{lemma}
There exists a constant $\beta_0>0$  such that
\begin{equation} \label{energia-domina}
\mathcal{E}_{\ell}(t) \geq  \beta_0 E_{\ell}(t) - Lm_F, \quad \forall \, t \geq  0.
\end{equation}
In addition, if $\ell \in (0,\pi/2L)$ then $\beta_0$ is independent of $\ell$. 
\end{lemma}

\proof We know from (\ref{energy-def}) and (\ref{new-F2}),
\begin{align*}
\mathcal{E}_{\ell}(t)
& \geq  E_{\ell}(t) - \beta \bigl( \Vert \varphi \Vert^2  + \Vert \psi \Vert^2 + \Vert w \Vert^2 \bigl) - L m_F \\
& \geq  \Big( 1 - \frac{2 \beta \gamma_3 L^2}{\pi^2}  \Big) E_{\ell}(t) - L m_F.
\end{align*}
Then from (\ref{Fc}) we take $\beta_0 = 1 - (2 \beta \gamma_3 L^2 / \pi^2)$. 
Finally, if $\ell \in (0, \pi/2L)$, then from Remark \ref{rem-ell} we can take $\gamma_3$ independent of $\ell$, and then $\beta_0$ 
is independent of $\ell$.  \qed


\subsection{Well-posedness} The existence of global weak and strong solutions to the Bresse system will be established through nonlinear semigroup theory. We shall write the system (\ref{P01})-(\ref{3ci}) as
an abstract Cauchy problem
\begin{equation} \label{PC}
\frac{d}{dt}y(t) -
\big( A_{\ell} + B \big) y(t) = \mathcal{F}(y(t)), \quad y(0)=y_{0},
\end{equation}
where
$$
y(t) =(\varphi(t), \psi(t), w(t),\tilde{\varphi}(t),\tilde{\psi}, \tilde{w}(t)) \in \mathcal{H}, \quad \tilde{\varphi} = \varphi_{t}, \; \tilde{\psi}=\psi_{t},\; \tilde{w}=w_t,
$$
and
$$
y_{0}=(\varphi_{0}, \psi_{0}, w_{0}, \varphi_{1},\psi_{1},w_{1}).
$$
To this end we take $A_{\ell}: D(A_{\ell}) \subset \mathcal{H} \to \mathcal{H}$,
$$
 A_{\ell}
    \left[ \begin{array}{c}   
      \varphi    \\
		  \psi\\
			w\\
      \tilde{\varphi} \\
      \tilde{\psi} \\
      \tilde{w}
    \end{array}\right] 
    =
  \left[ \begin{array}{c}
   \tilde{\varphi} \\
	   \tilde{\psi} \\
	   \tilde{w} \\
      \frac{k}{\rho_{1}}  (\varphi_{x}+\psi+\ell w)_{x} + \frac{k_{0}\ell}{\rho_{1}}(w_{x}-\ell\varphi)\\
      \frac{b}{\rho_{2}} \psi_{xx}- \frac{k}{\rho_{2}}(\varphi_{x} + \psi+\ell w) \\
      \frac{k_{0}}{\rho_{1}}(w_{x}-\ell\varphi)_{x} - \frac{k\ell}{\rho_{1}}(\varphi_{x}+\psi+\ell w) 
    \end{array}\right] ,
$$
with domain
$$
D(A_{\ell})=\bigl( H^{2}(0,L)\cap H^{1}_{0}(0,L)\bigl)^{3} \times H^{1}_{0}(0,L)^{3},
$$
and $B:\mathcal{H} \to \mathcal{H}$,
$$
B  \left[ \begin{array}{c}
      \varphi    \\
			\psi\\
			w\\
      \tilde{\varphi} \\
      \tilde{\psi} \\
      \tilde{w} \\
    \end{array}\right] =
  \left[ \begin{array}{c}
		0\\
		0\\
		0\\
      - g_{1}(\tilde{\varphi}) / \rho_1 \\
      - g_{2}(\tilde{\psi})/ \rho_{2} \\
     - g_{3}(\tilde{w})/\rho_{1} 
  \end{array}\right] ,
$$
with domain
$$
D(B)=\mathcal{H}.
$$
The forcing terms are represented by a nonlinear function
$\mathcal{F}: \mathcal{H} \to \mathcal{H}$
defined by
$$
\mathcal{F}  \left[ \begin{array}{c}
      \varphi    \\
			\psi\\
			w\\
      \tilde{\varphi} \\
      \tilde{\psi} \\
      \tilde{w} \\
    \end{array}\right]=
  \left[ \begin{array}{c}
  0\\ 
  0\\
  0\\ 
     -f_1(\varphi,\psi,w)/\rho_1  \\
      -f_2(\varphi,\psi,w)/\rho_2  \\
      -f_3(\varphi,\psi,w)/\rho_1 
 \end{array}\right].
$$

Our existence theorem is given in terms of the equivalent problem (\ref{PC}). 

\begin{theorem} [Well-posedness] \label{theo-existence} Assume that $\ell>0$ and the hypotheses
$(\ref{new-F0})$-$(\ref{HG2})$ hold.
Then for any initial data $y_{0} \in \mathcal{H}$ and $T>0$, problem $(\ref{PC})$ has 
a unique weak solution
$$
y \in C([0,T];\mathcal{H}), \;\; y(0)=y_0,
$$
given by
\begin{equation} \label{SolFormula}
y(t)= e^{t(A_{\ell}+B)}y_{0}+\int_{0}^{t}e^{(t-s)(A_{\ell}+B)}\mathcal{F}(y(s)) \normalfont{d} s, \;\; t \in [0,T], 
\end{equation}
and depends continuously on the initial data. 
In particular, if $y_0 \in D(A_{\ell})$ then the solution is strong.
\end{theorem}

\proof Under the hypotheses (\ref{HG1})-(\ref{HG2}) it was proved in \cite[Theorem 2.2]{charles-soriano} that $A_{\ell}+B$
is maximal monotone in $\mathcal{H}$. Then from standard theory the Cauchy
problem
\begin{equation} \label{CP-monotone}
\frac{d}{dt}y(t) - \big( A_{\ell} + B \big) y(t) = 0, \quad y(0)=y_{0},
\end{equation}
has a unique solution. We will show that system (\ref{PC}) is a locally Lipschitz perturbation of (\ref{CP-monotone}).
Then from classical results in \cite{barbu} (see a detailed proof in \cite[Theorem 7.2]{chueshov-eller-lasiecka}), we obtain a
local solution defined in an interval $[0,t_{\rm max})$ where, if $t_{max}< \infty$, then
\begin{equation} \label{limglobal}
\lim_{t \to \infty} \| y(t) \|_{\mathcal{H}} = + \infty .
\end{equation}


To show that operator $\mathcal{F}: \mathcal{H} \to \mathcal{H}$ is locally Lipschitz, let $G$ be a bounded set of $\mathcal{H}$ and
$y^1,y^2 \in G$. We can write $y^i = (z^i,z^i_t)$, where $z^j=(\varphi^i, \psi^i , w^i)$, $i=1,2$. 
Then from assumption (\ref{new-f1}) we obtain,
for $j=1,2,3$,
\begin{align*}
| f_j(z^1) -  f_j(z^2) | ^2
 = & \, |\nabla f_j ( \lambda z^1 + (1-\lambda) z^2)| ^2  \,  |z^1-z^2| ^2 \\
 \le & \, C_f ^2 \big(1
+ |\varphi^1|^{p-1}
+ |\varphi^2|^{p-1}
+ |\psi^1|^{p-1}
+ |\psi^2|^{p-1}
+ |w^1|^{p-1}
+ |w^2|^{p-1}
  \big)^2 \\ 
 & \, \times \big(
|\varphi^1 - \varphi^2|^2 +
|\psi^1 - \psi^2|^2 +
|w^1 - w^2|^2
\big) .  
\end{align*}
Then we infer that, for some $C_G >0$, 
$$
\int_0^L | f_j(z^1) -  f_j(z^2) |^2 dx
\le C_{G} \| z^1 - z^2\|_{(H^1_0)^3}^{2} \le C_{G} \| y^1 - y^2\|_{\mathcal{H}}^{2} .
$$
Summing this estimate on $j$ we obtain
$$
\| \mathcal{F}(y^1) - \mathcal{F}(y^2) \|_{\mathcal{H}}^2
\leq 3C_{G} \| y^1 - y^2 \|_{\mathcal{H}}^2 ,
$$
which shows that $\mathcal{F}$ is locally Lipschitz on $\mathcal{H}$. 

To see that the solution is global, that is, $t_{max}= \infty$, 
let $y(t)$ be a mild solution with initial data $y_0 \in D(A_{\ell}+B)$.
Then it is indeed a strong solution and so we can use energy inequality (\ref{energia-domina}) to conclude that
$$
 \| y(t) \|_{\mathcal{H}_{\ell}}^2 \leq \frac{2}{\beta_0}\left( \mathcal{E}_{\ell}(0) + L m_{f} \right), \;\; t \geq 0.
$$
By density, this inequality holds for mild solutions. Then clearly (\ref{limglobal}) does not hold and therefore
$t_{max} = \infty$.  

Finally, using (\ref{SolFormula}) we can check that for any initial 
data $y_0^1,y_0^2 \in \mathcal{H}$, the corresponding solutions 
$y^1,y^2$ satisfy 
$$
\| y^1(t)-y^2(t) \|_{\mathcal{H}}^2  \le C \| y_0^1-y_0^2 \|_{\mathcal{H}}^2, \quad 0<t<T, 
$$
which shows the continuous dependence on initial data.  \qed

\begin{remark} \label{rem-wp} \rm $(a)$ The existence of global solutions of 
(\ref{P01})-(\ref{3ci}) could be extended to the more general case involving damping terms with polynomial growth and no restriction near zero, as in \cite{peipei}. 
However, with respect to the existence of finite-dimensional global attractors,  
the method we use requires Lipschitz condition for the damping terms. 
Therefore, for brevity, we assumed conditions (\ref{HG1})-(\ref{HG2}) that were early  
considered in \cite{cavalcanti-timoshenko,charles-soriano}. 
$(b)$ The well-posedness of the Bresse system shows that 
its solution operator $S_{\ell}(t)$ is a $C^0$-semigroup on $\mathcal{H}$. Then we denote by 
$(\mathcal{H},S_{\ell}(t))$ the dynamical system generated 
by the problem $(\ref{P01})$-$(\ref{3ci})$. \qed 
\end{remark}

\setcounter{equation}{0}
\section{Singular limit}

Here we establish the Timoshenko limit of Bresse systems. 
With respect to the linear system $(\ref{BO1})$-$(\ref{BO3})$, if $\ell=0$, 
it produces the Timoshenko system $(\ref{TO1})$-$(\ref{TO2})$  
plus an independent wave equation in $w$. 
Therefore, in order to study the singular limit $\ell \to 0$ for the nonlinear model,  we shall need some compatibility condition.  
More precisely, we assume that 
$f_1,f_2$ do not depend on $w$, 
that is,  
\begin{equation} \label{singular50}
f_1(\varphi,\psi,w)=f_1(\varphi,\psi) \;\; \mbox{and} \;\; f_2(\varphi,\psi,w)=f_2(\varphi,\psi).
\end{equation}

\begin{remark} \label{rem-timo} \rm If the assumption  $(\ref{singular50})$ holds, then taking $\ell =0$, the same argument used in the proof of 
Theorem \ref{theo-existence} shows that 
Timoshenko system $(\ref{TS101})$-$(\ref{TS103})$ is well-posed 
in the phase space 
$$
\mathcal{H}_0 =  H^1_0(0,L)^2 \times L^2(0,L)^2.
$$
Its solution operator $S_0(t)$ generates a dynamical system denoted by $(\mathcal{H}_0,S_0(t))$. \qed
\end{remark}

\smallskip

\begin{theorem} \label{theo-singular} Assume that the  hypotheses of Theorem \ref{theo-existence} and  $(\ref{singular50})$ hold. 
Given any sequence $\{ \ell_{n} \}$ of positive numbers let 
$(\varphi^{n},\psi^{n}, w^{n})$ be the corresponding weak solution of the Bresse system $(\ref{P01})$-$(\ref{3ci})$, 
with $\ell = \ell_n$, and fixed initial data $(\varphi_0, \psi_0, w_{0}, \varphi_1,\psi_1, w_1)  \in  \mathcal{H}$.
Then if $\ell_{n} \to 0$ as $n \to \infty$, 
there exists $(\varphi,\psi,w ) $ such that, for any $T>0$, 
$$ 
(\varphi^{n},\psi^{n}, w^{n}) \stackrel{*}{\rightharpoonup} (\varphi,\psi,w ) \;\; \mbox{in} \;\; L^{\infty}(0,T; H^{1}_{0}(0,L)^{3}) ,
$$ 
$$ 
(\varphi^{n}_{t}, \psi^{n}_{t} , w^{n}_{t} )   \stackrel{*}{\rightharpoonup} 
( \varphi_{t},\psi_{t},w_{t} ) \;\; \mbox{in} \;\; L^{\infty}(0,T; L^{2}(0,L)^{3}),
$$ 
and $(\varphi,\psi)$ is a weak solution of the Timoshenko system $(\ref{TS101})$-$(\ref{TS103})$.
\end{theorem}

\proof The proof is divided into three parts. 

\medskip

\noindent $(i)$ A priori estimates: 
Since $\ell_n$ is uniformly bounded, there exists a positive constant $C_{0}$, 
such that,
$$
\mathcal{E}_{\ell_n}(0) = E_{\ell_n}(0) + 
\int_0^L F(\varphi_0,\psi_0,w_0)\, dx \le C_0, \quad \forall \, n.  
$$
Then, because $\mathcal{E}_{\ell_n}(t)$ is decreasing, we get from $(\ref{energia-domina})$,  
$$
E_{\ell_n}(t) \le \frac{1}{\beta_0} (C_0 + Lm_F), \quad \forall \, t \ge 0,
$$
From definition of $E_{\ell}(t)$ in (\ref{energy-linear}), 
we conclude that 
$$  
\{ \varphi^{n}_t \}, \;\; \{ \psi^{n}_t \}, \;\; \{ w^{n}_t \}, \;\; 
\{ \varphi^{n}_x + \psi^{n} + \ell_{n} w^{n} \}, \;\; \{ w^{n}_{x} - \ell_{n} \varphi^{n} \}, 
$$ 
are 
bounded in $L^{\infty}(0,T ; L^{2}(0,L))$ and $\{ \psi^{n} \}$ is bounded in $L^{\infty}(0,T; H^{1}_{0}(0,L))$. 
Let us show that $\{ \varphi^{n}\}, \{w^{n}\}$ are also bounded in $L^{\infty}(0,T; H^{1}_{0}(0,L))$. Indeed, from  
$$
\varphi^{n}(x,t) = \int_{0}^{t} \varphi^{n}_{t} (x,s) ds + \varphi_{0}(x),
$$
we infer that $\{\varphi^{n} \}$ is bounded in $L^{\infty}(0,T;  L^{2}(0,L))$. Now, using the relation
$$
w^{n}_{x} = (w^{n}_{x} - \ell_{n} \varphi^{n}) + \ell_{n} \varphi^{n}, 
$$
we find that  $\{w^{n} \}$ is bounded in $L^{\infty}(0,T; H^{1}_{0}(0,L))$. 
Similar arguments show that 
$\{ \varphi^{n} \}$ is bounded in $L^{\infty}(0,T; H^{1}_{0}(0,L))$. 

\medskip 

\noindent $(ii)$ Limits: Using a subsequence if necessary, there exist functions $\varphi,\psi, w, \vartheta_{1}, \vartheta_{2}$ such 
that
\begin{align}  \label{Conv1}
( \varphi^{n}, \psi^{n}, w^{n} ) \; \stackrel{*}{\rightharpoonup} ( \varphi, \psi, w )  \;\; \mbox{in} \;\; L^{\infty}(0,T; H^{1}_{0}(0,L)),
\end{align}
\begin{align}\label{Conv2}
( \varphi^{n}_t, \psi^{n}_t , w^{n}_t ) \; \stackrel{*}{\rightharpoonup} 
( \varphi_t, \psi_t, w_t ) \;\; \mbox{in} \;\; L^{\infty}(0,T ;L^{2}(0,L)),
\end{align}
\begin{align}\label{Conv3}
\{\varphi^{n}_x + \psi^{n} + \ell_{n} w^{n} \} \: \stackrel{*}{\rightharpoonup} \vartheta_{1}  \;\; \mbox{in} \;\; L^{\infty}(0,T; L^{2}(0,L)),
\end{align}
\begin{align}\label{Conv4}
\{ w^{n}_{x} - \ell_{n} \varphi^{n} \} 
\:\stackrel{*}{\rightharpoonup}  \vartheta_{2} \;\; \mbox{in} \;\; L^{\infty}(0,T; L^{2}(0,L)).
\end{align}
It follows from uniqueness of the weak limit, 
$$
\vartheta_{1} = \varphi_x + \psi \;\; \mbox{and} \;\; \vartheta_{2} = w_{x}.
$$
In addition, $(\ref{Conv1})$-$(\ref{Conv2})$ imply that  
\begin{align}\label{Conv5}
(\varphi^n,\psi^n) \to (\varphi,\psi) 
\;\; \mbox{in} \;\; L^{2}(0,T ;L^{2}(0,L)).
\end{align}
Now, from the definition of generalized solution for the Bresse system, we know that $( \varphi^{n} , \psi^{n}, w^{n} )$ 
satisfies 
\begin{align} \label{Weakform1} 
\rho_{1} \frac{d}{dt}  (\varphi_{t}^{n},\tilde{\varphi})
& + \rho_{2}\frac{d}{dt}(\psi^{n}_{t},\tilde{\psi})
+\rho_{1}\frac{d}{dt}(w^{n}_{t},\tilde{w})  
 + k\bigl((\varphi^{n}_{x}+\psi^{n}+\ell_{n} w^{n}),(\tilde{\varphi_{x}}  +\tilde{\psi}+\ell_{n}\tilde{w})\bigl) \nonumber \\ 
& + b(\psi^{n}_{x},\tilde{\psi}_{x}) 
  + k_{0}\bigl((w^{n}_{x}-\ell_{n}\varphi^{n}), (\tilde{w}_{x}-\ell_{n}\tilde{\varphi})\bigl) \nonumber \\ 
& + 
  N_{\varphi,\psi}^n + N_{w}^n
   =0, 
  \quad  \forall \, \tilde{\varphi},\tilde{\psi},\tilde{w} 
   \in H^{1}_{0}(0,L), 
 \end{align}  
where $N_{\varphi,\psi}^n$, 
$N_{w}^n$ denote nonlinear terms  
\begin{align*}
N_{\varphi,\psi}^n = 
\int_0^L f_1(\varphi^n,\psi^n)\tilde{\varphi} \,dx + \int_0^L f_2(\varphi^n,\psi^n)\tilde{\psi} \, dx 
+ \int_0^L g_1(\varphi_t^n)\tilde{\varphi} \,dx + \int_0^L g_2(\psi_t^n)\tilde{\psi} \, dx ,
\end{align*}
and 
\begin{align*}
N_{w}^n = \int_0^L f_3(\varphi^n, \psi^n, w^n)\tilde{w} \, dx + \int_0^L g_3(w_t^n)\tilde{w} \, dx .
\end{align*} 
Using $(\ref{Conv2})$ we infer that
$$
\int_0^L ( g_1(\varphi_t^n) - g_1(\varphi_t) )\tilde{\varphi} \, dx \to 0, \quad 
\int_0^L ( g_2(\psi_t^n) - g_2(\psi_t) )\tilde{\psi} \, dx \to 0. 
$$
Analogously, from (\ref{Conv5}) and (\ref{new-f1}) we infer that 
$$
\int_0^L ( f_1(\varphi^n,\psi^n) - f_1(\varphi,\psi) ) \tilde{\varphi} \, dx \to 0, \quad 
\int_0^L ( f_2(\varphi^n,\psi^n) - f_2(\varphi,\psi) ) \tilde{\psi} \, dx \to 0. 
$$
Then,  
taking test functions with $\tilde{w}=0$,    
we see that convergences (\ref{Conv1})-(\ref{Conv5}) imply that 
\begin{align}  \label{WeakTim}
\rho_{1}\frac{d}{dt}(\varphi_{t},\tilde{\varphi}) & +\rho_{2}\frac{d}{dt}(\psi_{t},\tilde{\psi})
+b(\psi_{x},\tilde{\psi}_{x})+k\bigl((\varphi_{x}+\psi),(\tilde{\varphi_{x}}+\tilde{\psi})\bigl) \nonumber \\
& + \int_0^L f_1(\varphi,\psi)\tilde{\varphi} \,dx + \int_0^L f_2(\varphi,\psi)\tilde{\psi} \, dx \nonumber \\ 
& + \int_0^L g_1(\varphi_t)\tilde{\varphi} \,dx + \int_0^L g_2(\psi_t)\tilde{\psi} \, dx 
=0, \quad  \forall \, \tilde{\varphi} ,\tilde{\psi} \in H^1_0(0,L).
\end{align}
This means that the limit $( \varphi, \psi)$ is a weak solution of the Timoshenko system (\ref{TS101})-(\ref{TS102}). 

\medskip

\noindent $(iii)$ Initial conditions: From (\ref{Conv1})-(\ref{Conv2}) we obtain (cf. \cite{simon}), 
\begin{align*}
( \varphi^{n}, \psi^{n} ) \rightarrow  ( \varphi, \psi ) \;\; \mbox{in} \;\; C([0,T],L^{2}(0,L)^{2}), 
\end{align*}
and therefore
\begin{equation} \label{Tip}
(\varphi(0), \psi(0) ) = (\varphi_0, \psi_0) .
\end{equation}
It remains to show that $(\varphi_t(0), \psi_t(0)) = (\varphi_1, \psi_1)$. 
To this end, we multiply $(\ref{Weakform1})$ 
by a test function 
$$
\theta \in H^{1}(0,T), \quad \theta(0)=1, \;\; \theta(T)=0,
$$  
and integrate over $[0,T]$. Taking also $\tilde{w}=0$, we find that
\begin{align*} 
-\rho_{1}(\varphi_{1},\tilde{\varphi}) & - \rho_{1}\int_{0}^{T}( \varphi_{t}^{n}, \tilde{\varphi}) \theta_{t} dt
- \rho_{2}(\psi_{1},\tilde{\psi}) - \rho_{2} \int_{0}^{T} (\psi_{t}^{n}, \tilde{\psi}) \theta_{t} dt + b \int_{0}^{T}(\psi^{n}_{x},\tilde{\psi}_{x}) \theta dt   \\
& + \, k \int_{0}^{T} \bigl( (\varphi^{n}_{x} + \psi^{n} + \ell_{n} w^{n}) , (\tilde{\varphi}_{x} + \tilde{\psi}) \bigl) \theta dt 
- k_{0} \int_{0}^{T} \ell_{n} (w_{x}^{n},\tilde{\varphi}) \theta dt\\
& +  k_{0} \int_{0}^{T} \ell_{n}^{2} (\varphi^{n},\tilde{\varphi}) \theta dt  
+ \int_0^T N_{\varphi, \psi}^n \, \theta \, dt =0,  \quad  \forall \, \tilde{\varphi} ,\tilde{\psi} \in H^1_0(0,L).
\end{align*}
Taking the limit $n \to  \infty$, we obtain 
\begin{align*} 
-\rho_{1}( & \varphi_{1}, \tilde{\varphi}) - \rho_{1} \int_{0}^{T}(\varphi_{t},\tilde{\varphi}) \theta_{t} dt
- \rho_{2} (\psi_{1},\tilde{\psi}) - \rho_{2} \int_{0}^{T}(\psi_{t},\tilde{\psi})\theta_{t} dt 
+ b\int_{0}^{T}(\psi_{x},\tilde{\psi}_{x}) \theta dt   \nonumber \\ 
& + k \int_{0}^{T} \big( (\varphi_{x} +\psi) , (\tilde{\varphi}_{x} + \tilde{\psi}) \big) \theta dt 
+ \int_0^T N^{\infty}_{\varphi,\psi} \, \theta \, dt = 0, \quad  \forall \, \tilde{\varphi} ,\tilde{\psi} \in H^1_0(0,L).
\end{align*}
On the other hand, multiplying (\ref{WeakTim}) by $\theta$ and integrating over $[0,T]$, we obtain 
\begin{align*} 
-\rho_{1}(& \varphi_{t} (0), \tilde{\varphi}) - \rho_{1} \int_{0}^{T}(\varphi_{t},\tilde{\varphi}) \theta_{t} dt
- \rho_{2} (\psi_{t}(0),\tilde{\psi}) - \rho_{2} \int_{0}^{T}(\psi_{t},\tilde{\psi})\theta_{t} dt 
+ b\int_{0}^{T}(\psi_{x},\tilde{\psi}_{x}) \theta dt  
 \nonumber \\ 
& + k \int_{0}^{T} \big( (\varphi_{x} +\psi) , (\tilde{\varphi}_{x} + \tilde{\psi}) \big) \theta dt 
+ \int_0^T N^{\infty}_{\varphi, \psi} \, \theta \, dt
= 0, \quad  \forall \, \tilde{\varphi} ,\tilde{\psi} \in H^1_0(0,L).
\end{align*}
The last two identities imply that 
\begin{equation} \label{Tiv}
(\varphi_{t}(0), \psi_{t}(0) ) = (\varphi_1, \psi_1 ).
\end{equation}

\noindent  Therefore (\ref{WeakTim}),(\ref{Tip}) and (\ref{Tiv}) show that the limit pair $(\varphi,\psi)$ is a solution of the Timoshenko system $(\ref{TS101})$-$(\ref{TS103})$. \qed 

\begin{remark} \rm We observe that the singular limit holds, as well as, for the linear problem with $f_i=0$, $g_i=0$, $i=1,2,3$. 
In this case, the energy is conservative and then $E_{\ell_n}(t)=E_{\ell_n}(0) \le C_0$, for all $t\ge 0$. \qed 
\end{remark}

\setcounter{equation}{0}
\section{Global attractors I} 

In this section we prove a first result on global attractors for Bresse systems. 
Some definitions and abstract results for global attractors are presented in the Appendix.

\begin{theorem} \label{theo-attractor1}
Under the hypotheses $(\ref{new-F0})$-$(\ref{HG2})$,
for each $\ell > 0$, the dynamical system $(\mathcal{H},S_{\ell}(t))$ generated by the problem
$(\ref{P01})$-$(\ref{3ci})$ has a global attractor $\mathbf{A}_{\ell}$. In addition, it is characterized by
$$
\mathbf{A}_{\ell}=\mathbb{M}_{+}(\mathcal{N}_{\ell}),
$$
where $\mathbb{M}_{+}(\mathcal{N}_{\ell})$ is the unstable manifold emanating from $\mathcal{N}_{\ell}$, 
the set of stationary points of $S_{\ell}(t)$.
\end{theorem}


The proof of this theorem is based on Theorem \ref{theo-A2}. We first show that the system is asymptotically compact. 

\begin{lemma} \label{lemma-ACI} Under the hypotheses 
Theorem \ref{theo-attractor1}, given a bounded subset $B$ of  $\mathcal{H}$,
let $S_{\ell}(t)y^{i} =(\varphi^{i},\psi^{i}, w^{i}, \varphi_{t}^{i}, \psi_{t}^{i}, w_{t}^{i})$ be,  with $i=1,2,$ two 
%
%
solutions of problem $(\ref{P01})$-$(\ref{3ci})$ with initial data $y^1,y^2 \in B$. Then, for every $\delta>0$, there exists a constant $C_{\delta, B}>0$ such that for $T>0$ sufficiently large, one has
\begin{equation} \label{quasistability1}
E_{\ell}(T) \leq \delta +\frac{C_{\delta, B}}{T}+C_{\delta, B}
\int_0^T \bigl( \|\varphi \|_{p+1} + \| \psi \|_{p+1} + \| w\|_{p+1} \bigl) dt,
\end{equation}
where
$\varphi=\varphi^{1}-\varphi^{2}$, $\psi=\psi^{1}-\psi^{2}$, $w=w^{1}-w^{2}$.
\end{lemma}

\proof For $u=u_1-u_2$ we use the following notation
$$
G_{i}(u)=g_{i}(u^{1})-g_{i}(u^{2})
\quad \mbox{and} \quad
F_i(u)=f_i(u^{1})-f_i(u^{2}), 
\quad i=1,2,3.
$$
Then $(\varphi, \psi, w,\varphi_{t}, \psi_{t}, w_{t})$ is the solution of the problem
\begin{align}
\rho_{1}\varphi_{tt} - k(\varphi_{x}+\psi+\ell w)_{x} - k_{0} \ell (w_{x}-\ell\varphi) & = - G_1(\varphi_t) - F_1(\varphi,\psi,w),   \label{B4} \\
\rho_{2}\psi_{tt} - b\psi_{xx}+k(\varphi_{x} + \psi+\ell w)
& = - G_2(\psi_t) - F_2(\varphi,\psi,w),   \label{B5} \\
\rho_{1}w_{tt} - k_{0}(w_{x}-\ell\varphi)_{x} + k\ell (\varphi_{x}+\psi+\ell w)
& = - G_3(w_t) - F_3(\varphi,\psi,w),  \label{B6}
\end{align}
with Dirichlet boundary condition and initial condition,
$$
( \varphi(0),\psi(0),w(0),\varphi_t(0),\psi_t(0),w_t(0) ) = y^1-y^2.
$$
Our objective is to obtain an estimate for $E_{\ell}(T)$. We begin by multiplying the equations (\ref{B4})-(\ref{B6})
by $\varphi$, $\psi$ and $w$, respectively, and integrate over $[0,T]\times[0,L]$. Then adding the kinetic energy, we obtain
\begin{align}\label{E12}
\int_{0}^{T}E_{\ell}(t)dt
 = & - \frac{1}{2}\int_{0}^{L} \bigl( \rho_{1}\varphi\varphi_t   +\rho_{2}\psi\psi_t + \rho_{1}w w_t \bigl) dx \Bigl|_{0}^{T} \nonumber \\
& + \int_{0}^{T}\int_{0}^{L}\bigl(\rho_{1}\varphi_t^{2} + \rho_{2}\psi_t^2 + \rho_{1}w_t^2\bigl) dxdt \nonumber \\
& - \frac{1}{2}\int_{0}^{T}\int_{0}^{L}\bigl(G_{1}(\varphi_t)\varphi + G_{2}(\psi_t)\psi + G_{3}(w_t)w\bigl) dxdt   \\
& - \frac{1}{2}\int_{0}^{T}\int_{0}^{L}\bigl( F_1(\varphi,\psi,w) \varphi + F_2(\varphi,\psi,w) \psi + F_3(\varphi,\psi,w) w\bigl) dx dt. \nonumber
\end{align}
We shall estimate the right-hand side of (\ref{E12}).

\noindent $(i)$ Boundary terms: Using H\"older's inequality and
norm inequality (\ref{norma-2}) there exists a constant $C>0$, independent of $\ell$,
such that
$$
\int_{0}^{L} \bigl( \rho_{1}\varphi \varphi_t + \rho_{2}\psi\psi_t+\rho_{1}w w_t \bigl) dx \leq C {E}_{\ell}(t), \quad \forall t \ge 0.
$$
Then we obtain
\begin{equation}\label{BT2}
\int_{0}^{L} \bigl( \rho_{1}\varphi \varphi_t + \rho_{2}\psi\psi_t+\rho_{1}w w_t \bigl) dx\Bigl|_{0}^{T}
\leq C \bigl( {E}_{\ell}(T)+{E}_{\ell}(0) \bigl).
\end{equation}

\noindent $(ii)$ Kinetic energy: Applying (\ref{HG4}), given $\delta > 0$, we have that there exists $C_{\delta}>0$ such that
$$
\int_{0}^{T}\int_{0}^{L}\varphi_t^{2}dx dt \leq TL\delta + C_{\delta}\int_{0}^{T}\int_{0}^{L}G_1(\varphi_t)\varphi_tdx dt.
$$
Same argument holds for $\iint \psi_t^{2}$ and $\iint w_t^{2}$ and therefore, given $\delta>0$, there exists $C_{\delta}>0$
such that
\begin{align} \label{kkk}
& \int_{0}^{T}\int_{0}^{L}\bigl( \rho_{1}\varphi_t^{2}+ \rho_{2}\psi_t^2 + \rho_{1}w_t^2\bigl) dx dt \nonumber \\
& \qquad \qquad  \leq  TL \delta + C_\delta \int_{0}^{T}\int_{0}^{L} \Bigl(G_1(\varphi_t) \varphi_t + G_2(\psi_t) \psi_t + G_3(w_t)w_t \Bigl) dxdt.
\end{align}

\noindent $(iii)$ Damping terms:  Let us consider the integral over $|\varphi_t| \leq 1$
and $|\varphi_t |>1$. Then (\ref{HG2}) implies 
that 
\begin{align*}
\int_{0}^{T} \! \int_{0}^{L} G_{1}(\varphi_t) \varphi dx dt
\leq & \, \int_{0}^{T} \! \int_{0}^{L}  \big( |g_1(\varphi_t^1)| + |g_1(\varphi_t^2)| \big) |\varphi| dx dt       \\
\leq & \, \int_{0}^{T} \! \int_{0}^{L} 2 \| g_1 \|_{L^{\infty}(-1,1)} |\varphi| dx dt  + \int_{0}^{T} \! \int_{0}^{L} M_1 \big( |\varphi_t^1| + |\varphi_t^2| \big) |\varphi| dx dt \\ 
\leq & \, C_{B} \int_0^T \| \varphi \|_{p+1} dt.
\end{align*}
The same argument holds for $\iint G_{2}(\psi_t)\psi$ and $\iint G_{3}(w_t)w$. Therefore we obtain the following estimate
\begin{align} \label{DT}
&  -\frac{1}{2} \int_{0}^{T}\int_{0}^{L}\Bigl(G_{1}(\varphi_t)\varphi + G_{2}(\psi_t)\psi + G_{3}(w_t)w\Bigl) dxdt \nonumber \\
& \qquad \qquad \leq
C_{B}\int_0^T \big( \| \varphi \|_{p+1} + \| \psi\|_{p+1}  + \| w \|_{p+1} \big) dt.
\end{align}
\noindent $(iv)$ Forcing terms: Using (\ref{new-f1}),
we have
\begin{align*}
\int_{0}^{T}\int_{0}^{L} F_1(\varphi,\psi,w) \varphi dx dt
& \leq  C \int_{0}^{T}\int_{0}^{L} C(\nabla f_{1}) (|\varphi|+|\psi|+|w|)|\varphi| dx dt \\
& \leq  C_{B} \int_0^T \bigl( \| \varphi \|_{p+1}^{2} +\| \psi \|_{p+1}^{2} +\| w \|_{p+1}^{2} \bigl)dt \\
& \leq  C_{B}  \int_0^T \bigl( \| \varphi \|_{p+1} +\| \psi \|_{p+1} + \| w \|_{p+1} \bigl) dt,
\end{align*}
where
$$
C( \nabla f_1 ) =C\bigl( 1 + |\varphi^1|^{p-1} + |\varphi^2|^{p-1}
+ |\psi^1 |^{p-1} + | \psi^2 |^{p-1}
+ |w^1|^{p-1} + | w^2 |^{p-1}\bigl).
$$
Analogous arguments with $\iint F_2 \psi$ and $\iint F_3 w$ imply that
\begin{align} \label{Fest1}
&& - \frac{1}{2}\int_{0}^{T}\int_{0}^{L}\Bigl( F_1(\varphi,\psi,w) \varphi
+ F_2(\varphi,\psi,w) \psi + F_3(\varphi,\psi,w) w\Bigl) dx dt \nonumber \\
&& \qquad \qquad \qquad \quad \leq C_{B} \int_0^T \bigl( \| \varphi \|_{p+1} + \| \psi\|_{p+1}  + \| w \|_{p+1} \bigl) dt.
\end{align}

\noindent $(v)$ An energy inequality: We multiply the equations (\ref{B4})-(\ref{B6}) by $\varphi_t$, $\psi_t$ and $w_t$, respectively,
and then integrate over $[s,T]\times [0,L]$. Then we find that
\begin{align} \label{EEnergy}
E_{\ell}(T)
 =  & \, E_{\ell}(s) -  \int_{s}^{T}\int_{0}^{L}
\bigl(G_{1}(\varphi_t)\varphi_{t}+G_{2}(\psi_t)\psi_{t} +G_{3}(w_t)w_{t} \bigl)  dxdt  \nonumber \\
 & - \int_{0}^{T}\int_{0}^{L}\bigl( F_1(\varphi,\psi,w) \varphi_t  + F_2(\varphi,\psi,w) \psi_t + F_3(\varphi,\psi,w) w_t \bigl) dx dt .
\end{align}
As before we see that
\begin{align*}
&& \int_{0}^{T}\int_{0}^{L}\bigl( F_1(\varphi,\psi,w) \varphi_t  + F_2(\varphi,\psi,w) \psi_t + F_3(\varphi,\psi,w) w_t \bigl) dx dt  \\
&& \qquad \qquad  \leq C_{B} \int_0^T \big( \| \varphi \|_{p+1} + \| \psi\|_{p+1}  + \| w \|_{p+1} \big) dt.
\end{align*}
Then identity (\ref{EEnergy}) gives
\begin{align} \label{GGG}
& \int_{0}^{T}\int_{0}^{L} \bigl( G_{1}(\varphi_t)\varphi_{t}+G_{2}(\psi_t)\psi_{t} +G_{3}(w_t)w_{t} \bigl) dxdt \nonumber \\
& \qquad \qquad \leq E_{\ell}(T) + E_{\ell}(0) + C_{B} \int_0^T \bigl( \| \varphi \|_{p+1} + \| \psi\|_{p+1}  + \| w \|_{p+1} \bigl) dt.
\end{align}

\noindent $(vi)$ Summarizing:  Inserting estimates (\ref{BT2})-(\ref{Fest1}) into (\ref{E12})
we have
\begin{align*}
\int_{0}^{T}E_{\ell}(t)dt
\leq &  \, C\bigl( E_{\ell}(T)+E_{\ell}(0)\bigl) + \delta \, T L  \\
& +  C_{\delta} \int_{0}^{T} \int_{0}^{L}\bigl(G_1(\varphi_t)\varphi_t+G_2(\psi_t)\psi_t +G_3(w_t)w_t \bigl) dx dt   \\
& + C_{B} \int_0^T \bigl( \| \varphi \|_{p+1} + \| \psi\|_{p+1} + \| w \|_{p+1} \bigl) dt.
\end{align*}
Then using (\ref{GGG}), and since $E(0) \leq C_{B}$, we can write
\begin{equation} \label{InE2}
\int_{0}^{T} E_{\ell}(t)dt
\leq \delta  T L + C_{\delta} E_{\ell}(T) +C_{\delta,B}
+ C_{\delta,B} \int_0^T \bigl( \| \varphi \|_{p+1} + \| \psi\|_{p+1}  + \| w \|_{p+1} \bigl) dt .
\end{equation}
Now integrating (\ref{EEnergy}) over $[0,T]$ with respect to the variable $s$ and
tanking into account that $G_1(\varphi_t)\varphi_t + G_2(\psi_t)\phi_t  + G_3(w_t)w_t \geq 0$, we have
\begin{equation} \label{temp}
TE_{\ell}(T)  \leq  \int_{0}^{T}{E}_{\ell}(t)dt + C_{B}T \int_0^T \bigl( \| \varphi \|_{p+1} + \| \psi\|_{p+1}  + \| w \|_{p+1} \bigl) dt.
\end{equation}
Combining (\ref{InE2}) and (\ref{temp}) we have
$$
T E_{\ell}(T)
\leq \delta T L + C_{\delta} E_{\ell}(T) + C_{\delta,B}+ C_{\delta,B} (1+T) \int_0^T \bigl(\| \varphi \|_{p+1} + \| \psi\|_{p+1}  + \| w \|_{p+1}\bigl) dt.
$$
Given $\delta >0$, we see that for $T$ sufficiently large (say $T> \max \{1 , \,   2C_{\delta} \}$) we can
write
$$
E_{\ell}(T) \leq   2 \delta L+ \frac{2C_{\delta,B}}{T}
+ 2 C_{\delta,B} \int_0^T \bigl( \| \varphi \|_{p+1} + \| \psi\|_{p+1}
+ \| w \|_{p+1} \bigl) dt .
$$
Then, renaming the constants we see that (\ref{quasistability1}) holds. \qed

\begin{lemma} \label{lemma-compact2}
Under the hypotheses of Theorem \ref{theo-attractor1} the system
$(\mathcal{H}, S_{\ell}(t))$ is asymptotically compact.
\end{lemma}

\proof Let $B$ be a positively invariant
bounded set of $\mathcal{H}_{\ell}$. Given $\varepsilon >0$,
we take $\delta$ sufficiently small and $T$ sufficiently large,
say,
$$
\delta < \frac{\varepsilon^2}{8} \quad \mbox{and} \quad \frac{C_{\delta,B}}{T} < \frac{\varepsilon^2}{8}. 
$$
Then from (\ref{quasistability1})
$$
\| S(T)y^1 - S(T)y^2 \|_{\mathcal{H}}
\le \varepsilon + \Phi_{\varepsilon,B,T} (y^1,y^2),
$$
where
$$
\Phi_{\varepsilon,B,T} (y^1,y^2) = 2 \sqrt{C_{\delta,B}}\Big( \int_0^T \big( \| \varphi^1 - \varphi^2 \|_{p+1} + \| \psi^1-\psi^2\|_{p+1}
+ \| w^1-w^2 \|_{p+1} \big) dt \Big)^{\frac{1}{2}} .
$$
Let us show that condition (\ref{contractive-f}) holds. Given $\{y^n\}$ in $B$, by positive invariance,
we see that $S_{\ell}(t)y^n = (\varphi^{n},\varphi_t^{n}, \psi^{n}, \psi_{t}^{n}, w^{i}, w_{t}^{n})$ is uniformly bounded in $\mathcal{H}$. Then
$$
(\varphi^n,\psi^n,w^n) \quad \mbox{is bounded in} \quad L^{\infty}\big(0,T;  H^1_0(0,L) ^3 \big),
$$
$$
(\varphi_t^n,\psi_t^n,w_t^n) \quad \mbox{is bounded in}
\quad L^{\infty}\big(0,T; L^2(0,L) ^3 \big),
$$
and therefore (\cite{simon}), 
$$
(\varphi^n,\psi^n,w^n) \quad \mbox{is pre-compact in} \quad C^0\big([0,T], L^{p+1}(0,L) ^3 \big).
$$
It follows that condition (\ref{contractive-f}) holds.
Then the asymptotic compactness follows from Theorem \ref{theo-A1}. \qed

\medskip

One of the assumptions of Theorem \ref{theo-A2} is that the set of stationary points is bounded.

\begin{lemma}
Under the hypotheses of Theorem \ref{theo-attractor1}, 
the set of equilibrium points $\mathcal{N}_{\ell}$ is bounded in $\mathcal{H}$.
\end{lemma}

\proof If $y \in \mathcal{N}_{\ell}$ we known that $y=(\varphi,\psi,w,0,0,0)$ and satisfies
\begin{align}
- k(\varphi_{x}+\psi+\ell w)_{x} - k_{0}\ell (w_{x}-\ell\varphi)  + f_1(\varphi, \psi,w) =0, & \label{ES1}\\
- b\psi_{xx}+k(\varphi_{x} + \psi+\ell w)+ f_2(\varphi, \psi,w)  =0 , & \label{ES2}\\
- k_{0}(w_{x}-\ell\varphi)_{x} + k\ell(\varphi_{x}+\psi+\ell w)  + f_3(\varphi, \psi,w) =0. & \label{ES3}
\end{align}
Multiplying (\ref{ES1}),(\ref{ES2}),(\ref{ES3}) by $\varphi, \psi, w$, respectively, and integrating over $[0,L]$, we obtain
\begin{align*} 
& \int_{0}^{L} \bigl( b\psi_{x}^{2} +k(\varphi_{x} + \psi+\ell w)^{2}+k_{0}(w_{x} + \ell\varphi)^{2}\bigl) dx \nonumber \\
& \qquad \qquad =  - \int_{0}^{L}\bigl( f_1(\varphi, \psi,w) \varphi +f_2(\varphi, \psi,w) \psi +f_2(\varphi, \psi,w) w\bigl)  dx.
\end{align*}
Then, using (\ref{norma-ok}), (\ref{new-fF}) and (\ref{new-F2}), we get
$$
\frac{1}{\gamma_3} \big(  \| \varphi_{x} \|^{2} + \|\psi_{x}\|^{2} + \| w_x \|^{2} \big) \leq  \frac{2\beta L^2}{\pi^2}  \big(  \| \varphi_{x} \|^{2} + \|\psi_{x}\|^{2} + \| w_x \|^{2} \big) 
+  2m_FL , 
$$
and consequently,  
$$
\Big( 1 - \frac{2\beta L^2\gamma_3}{\pi^2}  \Big) \big(  \| \varphi_{x} \|^{2} + \|\psi_{x}\|^{2} + \| w_x \|^{2} \big) \leq   2m_FL\gamma_3. 
$$
Therefore the set $\mathcal{N}_{\ell}$ is bounded in $\mathcal{H}$. \qed

\paragraph{Proof of Theorem \ref{theo-attractor1} (completion).}  
We already known that the system is asymptotically compact and the set of its stationary points 
$\mathcal{N}_{\ell}$ is bounded. To apply Theorem \ref{theo-A2}, it remains to 
show that the dynamical system $(\mathcal{H},S_{\ell}(t))$ is gradient and satisfies condition 
(\ref{SSE}). Indeed, we can take the energy functional $\mathcal{E}_{\ell}$ as 
a Lyapunov function $\Phi$, since $t \rightarrow \Phi(S_{\ell}(t)y)$ 
is then strictly decreasing for any $y \in \mathcal{H}$. 
Moreover, from (\ref{energy-def}) and (\ref{new-f1}) we see that 
$\mathcal{E}_{\ell}(t) \le \Vert y(t) \Vert_{\mathcal{H}_{\ell}}^2 
+ C (1+ \Vert y(t) \Vert_{\mathcal{H}_{\ell}}^{p+1})$. Then 
$$
\Phi(y) \to \infty \;\; \mbox{implies that} \;\; \Vert y \Vert_{\mathcal{H}_{\ell}} \to \infty .   
$$
On the other hand, the inequality $(\ref{energia-domina})$ implies that 
$E_{\ell}(t) \le \frac{1}{\beta_0} (\mathcal{E}_{\ell} (t) + Lm_F )$, 
and then
$$
\Vert y \Vert_{\mathcal{H}_{\ell}} \to \infty \;\; \mbox{implies that} \;\;  \Phi(y) \to \infty.  
$$
Therefore all the assumptions of Theorem \ref{theo-A2} are fulfilled and consequently 
the system $(\mathcal{H},S_{\ell}(t))$ has a global attractor $\mathbf{A}_{\ell} =\mathbb{M}_{+}(\mathcal{N}_{\ell})$. 
This ends the proof of Theorem \ref{theo-attractor1}. \qed 

\begin{remark} 
\rm The existence of a global attractor implies 
that the system has a bounded absorbing set. 
But in principle it depends on $\ell$. We shall construct an absorbing set which is uniformly bounded 
for $\ell \in [0,\ell_0]$, with $\ell_0$ small. This will be used in Theorem \ref{theo-upper}. \qed
\end{remark}

\begin{lemma} \label{lemma-AbsobSet} 
Under hypotheses of Theorem \ref{theo-attractor1},   
with $\ell \in (0, \pi/2L)$, the system $(\mathcal{H},S_{\ell}(t))$ has a bounded 
absorbing set $\mathcal{B}$ independent of $\ell$. 
\end{lemma}

\proof Multiply the equations (\ref{P01})-(\ref{P03}) by 
$\varphi$, $\psi$ and $w$, respectively, and integrate over $[0,L]\times [0,T]$. We obtain 
\begin{align} \label{E1}
\int_{0}^{T}& \Big(b\|\psi_x\|^2 + k\|\varphi_x + \psi+\ell w\|^2+ k_0\|w_x-\ell\varphi\|^2  \Big) dt 
 + \int_0^T \!\! \int_{0}^{L} \nabla F(\varphi,\psi,w) 
\!\cdot\! (\varphi,\psi,w) dx dt \nonumber \\ 
& = - \int_{0}^{L}\bigl(\rho_{1}\varphi\varphi_t+\rho_{2}\psi\psi_t+\rho_{1}w w_t\bigl) dx \Bigl|_{0}^{T}
+\int_{0}^{T}\int_{0}^{L}\bigl(\rho_{1}\varphi_t^{2}+\rho_{2}\psi_t^2+\rho_{1}w_t^2\bigl)dx dt  \\
&\quad \: -\int_{0}^{T}\int_{0}^{L}\bigl(g_{1}(\varphi_t)\varphi+g_{2}(\psi_t)\psi + g_{3}(w_t)w\bigl) dx dt. \nonumber 
\end{align}
Inequality (\ref{new-fF}) together with (\ref{norma-ok}) and (\ref{Fc}) implies that  
\begin{align*}
& \int_0^T  \int_{0}^{L}  \nabla F(\varphi,\psi,w)\! \cdot \! (\varphi,\psi,w) dx dt  \\
& \qquad \geq  \int_0^T \int_{0}^{L} F(\varphi,\psi,w) dx dt - \int_0^T \int_{0}^{L} \beta \bigl( |\varphi|^2  + |\psi |^2 + | w |^2 \bigl)dx dt - TLm_F \\
& \qquad \geq  \int_0^T \int_{0}^{L} F(\varphi,\psi,w) dx dt  
- \frac{1}{2}\int_0^T \bigl( b\|\psi_x\|^2 + k\|\varphi_x + \psi+\ell w\|^2+ k_0\|w_x-\ell\varphi\|^2  \bigl)dt \\ 
& \qquad \quad \, - TLm_F.
\end{align*}
Inserting this inequality into (\ref{E1}) and adding the kinetic energy, we infer that
\begin{align}
\label{E2}
\int_{0}^{T}\mathcal{E}_{\ell}(t) dt \leq 
& - \int_{0}^{L}\bigl(\rho_{1}\varphi\varphi_t + \rho_{2}\psi\psi_t+\rho_{1}w w_t\bigl) dx \Bigl|_{0}^{T} 
+ \frac{3}{2}  \int_{0}^{T}\int_{0}^{L}\bigl(\rho_{1}\varphi_t^{2}+\rho_{2}\psi_t^2+\rho_{1}w_t^2\bigl) dx dt \nonumber \\
& - \int_{0}^{T}\int_{0}^{L}\bigl(g_{1}(\varphi_t)\varphi + g_{2}(\psi_t)\psi+g_{3}(w_t)w\bigl)dx dt + TLm_F.  
\end{align}

\noindent In the following we will estimate the terms on the right-hand side of (\ref{E2}). 

\medskip

\noindent {$(i)$} Estimates for the boundary terms:  Young's inequality and (\ref{norma-ok}) imply 
that  
\begin{align*}
-\int_{0}^{L} \bigl( \rho_{1}\varphi \varphi_t + \rho_{2}\psi\psi_t+\rho_{1}w w_t \bigl) dx \leq  C E_{\ell}(t)
\end{align*}
for some constant $C>0$, independent of $T$ and $\ell$. 
Using inequality (\ref{energia-domina}) we obtain 
\begin{align*}
-\int_{0}^{L} \bigl( \rho_{1}\varphi \varphi_t + \rho_{2}\psi\psi_t+\rho_{1}w w_t \bigl)dx \leq   \frac{2}{\beta_0} \mathcal{E}_{\ell}(t)  + \frac{2Lm_F}{\beta_0} .
\end{align*}
Noting that $\beta_0$ does not depend on $\ell$, there exists $C_{1}>0$, independent of $T$ and $\ell$, such that 
\begin{align}\label{BT}
-\int_{0}^{L} \bigl(\rho_{1}\varphi \varphi_t + \rho_{2}\psi\psi_t+\rho_{1}w w_t \bigl) dx \Bigl|_{0}^{T}
\leq  C_1 \bigl( \mathcal{E}_{\ell}(T)+\mathcal{E}_{\ell}(0)\bigl) + C_1.
\end{align}

\medskip

\noindent {$(ii)$} Estimates for the damping terms:  Using Young's inequality  we have that
\begin{align*}
& - \int_{0}^{T}\int_{0}^{L}\bigl(g_{1}(\varphi_t)\varphi+g_{2}(\psi_t)\psi+g_{3}(w_t)w\bigl) dx dt \\
& \qquad \leq  \frac{1}{2} \int_0^T E_{\ell}(t)
dt + C \int_{0}^{T}\int_{0}^{L}\bigl(g_{1}(\varphi_t)^2 +g_{2}(\psi_t)^2 + g_{3}(w_t)^2 \bigl) dx dt. 
\end{align*}
From assumption (\ref{HG2}), for $i=1,2,3$, we get that
$$
\int_{0}^{T}\int_{|u|\leq  1} g_{i}(u)^2 dx dt \leq  \max \{g(-1)^2 , g(1)^2 \}LT
$$
and
$$
\int_{0}^{T}\int_{|u| > 1} g_{i}(u)^2 dx dt  \leq  M_i \int_{0}^{T}\int_{0}^{L} g_i(u)u dx dt.
$$
Then there exists a constant $C_2> 0$, independent of $T$ and $\ell$,  such that 
\begin{align} \label{DD}
& - \int_{0}^{T}\int_{0}^{L}\bigl(g_{1}(\varphi_t)\varphi+g_{2}(\psi_t)\psi+g_{3}(w_t)w\bigl) dx dt  \nonumber \\
& \qquad \leq \frac{1}{2} \int_{0}^{T}\mathcal{E}_{\ell}(t) dt + C_2 \int_{0}^{T}\int_{0}^{L}\bigl(g_{1}(\varphi_t)\varphi_t + g_{2}(\psi_t)\psi_t + g_{2}(w_t)w_t \bigl) dx dt  
+ C_2 T.
\end{align}

\medskip

\noindent {$(iii)$} Estimates for the kinetic energy:  Firstly we note that using (\ref{HG2}) we have
\begin{align*}
\int_{0}^{T}\int_{|\varphi_t| >1} \varphi_t^{2} dxdt \leq  \frac{1}{m_1} \int_{0}^{T}\int_{0}^{L} g_{1}(\varphi_t)\varphi_t dxdt 
\end{align*}
and then 
\begin{align*}
\int_{0}^{T}\int_{0}^{L}\varphi_t^{2} dxdt  &=\int_{0}^{T}\int_{|\varphi_t| >1}\varphi_t^{2}
dxdt + \int_{0}^{T}\int_{|\varphi_t| \leq 1}\varphi_t^{2}dx dt\\
&\leq \frac{1}{m_1} \int_{0}^{T}\int_{0}^{L} g_{1}(\varphi_t)\varphi_t dx dt + TL.
\end{align*}
Similar estimate holds for $\iint \psi_t^2$ and $\iint w_t^2$. Therefore there exists a  constant $C_3>0$, 
independent of $T$ and $\ell$, 
such that 
\begin{align} \label{KE}
\begin{split}
& \frac{3}{2}  \int_{0}^{T}\int_{0}^{L}\bigl(\rho_{1}\varphi_t^{2}+\rho_{2}\psi_t^2+\rho_{1}w_t^2\bigl) dx dt \\ 
& \qquad \qquad \leq  C_3 \int_{0}^{T}\int_{0}^{L}\bigl(g_1(\varphi_t) \varphi_t + g_2(\psi_t) \psi_t  + g_3(w_t) w_t  \bigl) dx dt + C_3 T.
\end{split}
\end{align}

\medskip

\noindent {$(iv)$} Inserting estimates (\ref{BT}), (\ref{DD}), (\ref{KE}) into (\ref{E2}) we obtain
\begin{align*}
\frac{1}{2} \int_{0}^{T} \mathcal{E}_{\ell}(t) dt \leq 
& \; C_1 \bigl( \mathcal{E}_{\ell}(T) + \mathcal{E}_{\ell}(0)\bigl) \\ 
& + (C_2+C_3) \int_{0}^{T}\int_{0}^{L}\bigl(g_{1}(\varphi_t)\varphi_t+g_{1}(\psi_t)\psi_t+g_{1}(w_t)w_t \bigl) dx dt \\ 
& + C_{1} + (C_2+C_3) T.
\end{align*}
Using the energy identity (\ref{energy-ID}) and noting that $\mathcal{E}(T) \leq \mathcal{E}(t)$ in the left-hand side 
integral,  
\begin{align*}
\frac{T}{2}\mathcal{E}_{\ell}(T)  
\leq \big( C_{1} - C_2 - C_3 \big) \mathcal{E}_{\ell}(T) + \big( C_{1}+ C_2 + C_{3}) \mathcal{E}_{\ell}(0)+ C_1 + (C_2+C_3)T .
\end{align*}
Taking $T$ sufficiently ($T> 2C_1$) we can write
\begin{align} \label{absorbing-quasi}
\mathcal{E}_{\ell} (T) \leq \gamma_T \mathcal{E}_{\ell}(0) + K_T, 
\end{align}
where 
$$
\gamma_T = \frac{2(C_1+C_2+C_3)}{T-2(C_1-C_2-C_3)} <1 , \quad  K_T=\frac{2C_1+2(C_2+C_3)T }{T- 2(C_1-C_2-C_3)}>0 .
$$
From (\ref{absorbing-quasi}) and well-known argument shows that  
\begin{equation} \label{423}
\mathcal{E}_{\ell}(t) \leq \gamma \mathcal{E}_{\ell}(0)e^{-\alpha t} + \frac{K_T}{1-\gamma_T} , \quad \forall \, t \ge 0,   
\end{equation}
for some $\alpha, \gamma >0$. For completeness we sketch its proof here.  Indeed, 
the same argument can be repeated for any interval $[mT,(m+1)T]$, $m\in \mathbb{N}$. 
Then 
\begin{align*}
\mathcal{E}_{\ell}(mT) 
& \leq \gamma_{T}  \mathcal{E}_{\ell}((m-1)T) + K_T \\
& \leq \gamma_T^{m} \mathcal{E}_{\ell}(0) + \Big(  \sum_{j=0}^{m-1} \gamma_T^{j}  \Big) K_T \\
& \leq \gamma_T^{m} \mathcal{E}_{\ell}(0) + \frac{K_T}{1-\gamma_T} . \quad (\mbox{since $\gamma_T < 1$})
\end{align*}
Now given $t \geqslant 0$, there exits $m \in \mathbb{N}$ and $r\in [0,T)$ such that $t = m T + r$. Then
$$
\mathcal{E}_{\ell}(t)  \leq \mathcal{E}_{\ell}(mT) \leq \gamma_T^m \mathcal{E}_{\ell}(0)+\frac{K_T}{1-\gamma_T}  .
$$
It follows that 
$$
\mathcal{E}_{\ell}(t)  \leq  \gamma_T^{\frac{t-r}{T}} \mathcal{E}_{\ell}(0)+\frac{K_T}{1-\gamma_T}  
 \leq  \gamma_T^{-1} \gamma_T^{\frac{t}{T}} \mathcal{E}_{\ell}(0)+\frac{K_T}{1-\gamma_T} .
$$
Therefore choosing $\gamma = \gamma_T^{-1}$ and $\alpha = - \ln (\gamma_T)/T$ we obtain (\ref{423}). 

\medskip 

\noindent {$(v)$} Conclusion: We observe that combining (\ref{423}) and (\ref{energia-domina}) yields 
$$
\| S(t) y_0 \|_{\mathcal{H}_{\ell}}^2 \le \frac{2 \gamma}{\beta_0} \mathcal{E}_{\ell}(0) e^{- \alpha t} + \frac{2Lm_F K_T}{\beta_0(1-\gamma_T)},  
$$
and then clearly any closed ball $\overline{B}_{\mathcal{H}}(0,R_0)$ with $R_0^2 > \frac{2Lm_F K_T}{\beta_0(1-\gamma_T)}$ 
is a bounded absorbing set, not depending 
on $\ell$. \qed

\setcounter{equation}{0}
\section{Global attractors II} 

In this section we assume that damping terms satisfy condition (\ref{HG3}). Then we show that the global attractor obtained in Theorem \ref{theo-attractor1} 
has further properties. 

\begin{theorem}\label{theo-attractor2} Under the hypotheses of Theorem \ref{theo-attractor1}, with $(\ref{HG2})$ replaced by $(\ref{HG3})$, one has: 

\noindent $(i)$ The global attractor $\mathbf{A}_{\ell}$ has finite fractal dimension.

\medskip

\noindent $(ii)$ Any full trajectory $(\varphi(t),  \psi(t),  w(t),\varphi_t(t),\psi_t(t) , w_t(t) )$  
inside the attractor $\mathbf{A}_{\ell}$, has further regularity
\begin{equation} \label{reg100}
\Vert (\varphi, \psi, w)\Vert_{(H^{2})^3} + \Vert (\varphi_t, \psi_t, w_t) \Vert_{(H^{1}_0)^3}  + \Vert (\varphi_{tt}, \psi_{tt}, w_{tt}) \Vert_{(L^2)^3}  \le C_{\ell},  
\end{equation}
for some $C_{\ell}>0$. 

\medskip

\noindent $(iii)$ The dynamical system $(\mathcal{H},S_{\ell}(t))$
possesses a generalized exponential attractor $\mathbf{A}_{\ell}^{\rm exp}$,
with finite fractal dimension in a extended space 
$\widetilde{\mathcal{H}}_{-\eta}$, defined as interpolation of 
$$
\widetilde{\mathcal{H}}_0 := \mathcal{H} \quad  \mbox{and} \quad \widetilde{\mathcal{H}}_{-1}:= L^2(0,L)^3 \times H^{-1}(0,L)^3,
$$
for any $\eta \in (0,1]$.
\end{theorem}

\begin{remark}  
\rm As discussed in Remark \ref{rem-timo}, we can prove an analogous result for the Timoshenko system,  
that is, the dynamical system $(\mathcal{H}_0,S_ 0(t))$ generated by $(\ref{TS101})$-$(\ref{TS103})$ 
has a regular global attractor $\mathbf{A}_0$ in $\mathcal{H}_0 = H^1_0(0,L)^2\times L^2(0,L)^2$,  
with finite fractal dimension. \qed 
\end{remark}

The proof of Theorem \ref{theo-attractor2} relies on the properties of quasi-stable systems.

\subsection{Quasistability}

\begin{lemma} \label{lemma-QEI} In the context of Lemma \ref{lemma-ACI}, with  $(\ref{HG2})$ replaced by $(\ref{HG3})$, 
given a bounded invariant set $B$, there exist  constants $\alpha_{B},  \gamma_{B}, C_{B}>0$, such
that
\begin{equation} \label{q-stable}
E_{\ell}(t) \leq \gamma_{B}  E_{\ell}(0)e^{-\alpha_{B} t} + C_{B} \sup_{\sigma \in [0,t]} \big( \|\varphi(\sigma) \|_{2p}^{2} + \|\psi(\sigma) \|_{2p}^{2} 
+ \| w(\sigma) \|_{2p}^{2} \bigl).
\end{equation}
\end{lemma}

\proof We begin as in the proof of Lemma \ref{lemma-ACI}, since $(\ref{HG3})$ implies $(\ref{HG2})$. Then we only need estimate
the right-hand side of (\ref{E12}). Here $C>0$ will represent several constants independent of $B$ or $t$.  

\medskip

\noindent $(i)$ First remarks:  We observe that estimate (\ref{BT2}) holds unchanged. We also observe that (\ref{Fest1}) can be changed to
\begin{align} \label{Fest2}
&& - \frac{1}{2}\int_{0}^{T}\int_{0}^{L}\bigl(F_1(\varphi,\psi,w) \varphi
+ F_2(\varphi,\psi,w) \psi + F_3(\varphi,\psi,w) w\bigl)dx dt \nonumber \\
&& \qquad \qquad \qquad  \leq  C_{B} \int_0^T \bigl( \| \varphi \|_{2p}^2
+ \| \psi\|_{2p}^2  + \| w \|_{2p}^2 \bigl) dt.
\end{align}

\noindent $(ii)$ Role of (\ref{HG3}): Now since (\ref{HG5}) holds  we see that estimate (\ref{kkk})
becomes
$$
 \int_{0}^{T} \!\! \int_{0}^{L}\bigl( \rho_{1}\varphi_t^{2}+ \rho_{2}\psi_t^2 + \rho_{1}w_t^2\bigl) dx dt
 \leq C \int_{0}^{T} \!\! \int_{0}^{L} \bigl(G_1(\varphi_t) \varphi_t + G_2(\psi_t) \psi_t + G_3(w_t)w_t \bigl)  dx dt .
$$
In addition, (\ref{HG3}) implies that $| g_i(u)-g_i(v) | \le M_i |u-v|$ for all
$u,v \in \mathbb{R}$. Then
$$
\int_{0}^{T}\int_{0}^{L} G_{1}(\varphi_t)\varphi \, dx dt
\leq  \frac{1}{6}\int_0^T \| \varphi_t \|_{2}^2 \, dt dx
+ C\int_0^T \| \varphi \|_{2}^2 \,dt.
$$
Applying
the same argument to $\iint G_2(\psi_t)\psi$ and $\iint G_3(w_t)w$  we
infer that (\ref{DT}) becomes
\begin{align} \label{DT2}
& -\frac{1}{2} \int_{0}^{T}\int_{0}^{L}\bigl(G_{1}(\varphi_t)\varphi + G_{2}(\psi_t)\psi + G_{3}(w_t)w\bigl) dxdt \nonumber  \\
& \qquad \qquad \leq  \frac{1}{2} \int_{0}^{T} E_{\ell}(t)dt +
C \int_0^T \bigl( \| \varphi \|_{2p}^2 + \| \psi\|_{2p}^2
+ \| w \|_{2p}^2  \bigl) dt.
\end{align}

\noindent $(iii)$ First energy inequality: Using the inequalities (\ref{q-stable})-(\ref{DT2}) into (\ref{E12}), we obtain that
\begin{align} \label{InEE}
\int_{0}^{T}E_{\ell}(t)dt
\leq & \, C \bigl[E_{\ell}(T)+E_{\ell}(0)\bigl] + C_{B} \int_{0}^{T}\int_{0}^{L}\bigl(G_1(\varphi_t)\varphi_t+ G_2(\psi_t)\psi_t + G_3(w_t)w_t \bigl) dx dt \nonumber \\
& + C_{B} \int_0^T \bigl( \| \varphi \|_{2p}^2 + \| \psi\|_{2p}^2 + \| w \|_{2p}^2  \bigl) dt.
\end{align}

\noindent $(iv)$ Damping estimate: The energy identity (\ref{EEnergy})  implies that
%
%
%
\begin{align} \label{InEE333}
& \int_{0}^{T}\int_{0}^{L}\bigl(G_1(\varphi_t)\varphi_t + G_2(\psi_t)\psi_t +G_3(w_t)w_t \bigl)dxdt \nonumber \\
& \qquad  =  - \int_{0}^{T}\int_{0}^{L} \bigl( F_1(\varphi,\psi,w) \varphi_t 
+ F_2(\varphi,\psi,w) \psi_t + F_3(\varphi,\psi,w) w_t \bigl) dxdt. \nonumber \\
& \qquad \quad  + \, E_{\ell}(0) - E_{\ell}(T) .
\end{align}
Let 
us estimate the forcing terms. Note that, for $\epsilon>0$, 
\begin{align*}
 \int_{0}^{L} F_1(\varphi,\psi,w) \varphi_t  dx
 &\leq 
C(\nabla f_{1})
( \| \varphi \|_{2p} + \| \psi \|_{2p}  +\| w \|_{2p}) \| \varphi_t \|_{2}  \\
& \leq 
 \frac{\epsilon}{3T}\| \varphi_t \|_2^2 + T \:C_{\epsilon,B}
\big( \| \varphi \|_{2p}^2 + \| \psi\|_{2p}^2
+ \| w \|_{2p}^2  \big).
\end{align*}
where
$$
C(\nabla f_{1})=
C\big(1 + \| \varphi^{1} \|_{2p}^{p-1} +\| \varphi^{2} \|_{2p}^{p-1}
+  \| \psi^{1} \|_{2p}^{p-1} +\| \psi^{2} \|_{2p}^{p-1}+ \| w^{1} \|_{2p}^{p-1}+ \| w^{2} \|_{2p}^{p-1} \big) .
$$
Similar estimate holds for $\int F_2(\varphi,\psi,w) \psi_t$ and $\int F_3(\varphi,\psi,w) w_t$. 
Then we obtain
\begin{align*}
& \int_{0}^{L}
 \bigl( F_1(\varphi,\psi,w) \varphi_t
 + F_2(\varphi,\psi,w) \psi_t + F_3(\varphi,\psi,w) w_t \bigl) dx \\
& \qquad \qquad \leq
 \frac{\epsilon}{T} \big( \| \varphi_t \|_2^2 + \| \psi_t \|_2^2
+ \| w_t \|_2^2 \big) + T \:C_{\epsilon,B}
\big( \| \varphi \|_{2p}^2 + \| \psi\|_{2p}^2
+ \| w \|_{2p}^2  \big)
\end{align*}

\noindent and we can write
\begin{align}\label{FEstk}
& \int_{0}^{T}\int_{s}^{T}\int_{0}^{L}
 \bigl( F_1(\varphi,\psi,w) \varphi_t
 + F_2(\varphi,\psi,w) \psi_t + F_3(\varphi,\psi,w) w_t \bigl) dxdtds
 \nonumber \\
& \qquad \qquad \leq \epsilon\int_{0}^{T}E_{\ell}(t)dt +
C_{\epsilon, B, T} \int_0^T \bigl( \| \varphi \|_{2p}^2 + \| \psi\|_{2p}^2
+ \| w \|_{2p}^2 \bigl) dt.
\end{align}

\noindent This inequality together with (\ref{InEE333}) results that
\begin{align} \label{DampEstK}
& \int_{0}^{T}\int_{0}^{L} \bigl(G_1(\varphi_t)\varphi_t
 + G_2(\psi_t)\psi_t + G_3(w_t)w_t \bigl) dx dt  \nonumber \\
& \qquad \qquad \leq E_{\ell}(0) - E_{\ell}(T)   + \epsilon \int_0^T E_{\ell}(t)dt   \nonumber \\ 
& \qquad \qquad \quad \, 
+ C_{\epsilon, B, T} \int_0^T \bigl( \| \varphi \|_{2p}^2 + \| \psi\|_{2p}^2
+ \| w \|_{2p}^2  \bigl) dt. 
\end{align}

\noindent $(v)$ Second energy inequality:  Applying the damping estimate $(\ref{DampEstK})$ in $(\ref{InEE} )$ we obtain, for $\epsilon$ small enough,
\begin{equation} \label{resumo}
\int_{0}^{T}E_{\ell}(t)dt \leq  (C-C_{B}) E_{\ell}(T)+(C+C_{B})E_{\ell}(0)
+C_{B} \int_0^T \bigl( \| \varphi \|_{2p}^2 + \| \psi\|_{2p}^2 + \| w \|_{2p}^2  \bigl) dt.
\end{equation}

\noindent $(vi)$ Estimating $E(T)$: Integrating the energy identity (\ref{EEnergy}) it follows
that
\begin{align*} 
T E_{\ell}(T)
 = & \int_{0}^{T}E_{\ell}(t)dt
- \int_{0}^{T}\int_{s}^{T}\int_{0}^{L}\bigl(G_1(\varphi_t)\varphi_t +G_2(\psi_t)\psi_t +G_3(w_t)w_t \bigl)dxdtds  \\
& - \, \int_{0}^{T}\int_{s}^{T}\int_{0}^{L} \bigl( F_1(\varphi,\psi,w) \varphi_t + F_2(\varphi,\psi,w) \psi_t + F_3(\varphi,\psi,w) w_t \bigl) dxdtds.
\end{align*}
Taking 
into account that $G_1(\varphi_t)\varphi_t + G_2(\psi_t)\phi_t  + G_3(w_t)w_t \geq 0$ and using the estimate $(\ref{FEstk})$ we obtain
\begin{equation} \label{A21}
T E_{\ell}(T)
\leq 2 \int_{0}^{T}E_{\ell}(t)dt
+ C_{B,T} \int_0^T \bigl( \| \varphi \|_{2p}^2 + \| \psi\|_{2p}^2 + \| w \|_{2p}^2 \bigl) dt.
\end{equation}

\noindent $(vii)$ Concluding: Inserting (\ref{resumo}) into (\ref{A21})
we obtain
$$
T E_{\ell}(T) \leq 2(C-C_{B}) E_{\ell}(T)+ 2(C+C_{B})E_{\ell}(0)+ C_{B,T}\int_0^T \bigl( \| \varphi \|_{2p}^2 + \| \psi\|_{2p}^2
+ \| w \|_{2p}^2 \bigl) dt.
$$
Taking $T>  4C$ we can write
$$
E_{\ell}(T) \leq \gamma_T E_{\ell}(0)+ C_{B,T} \sup_{\sigma \in [0,T]}\big( \| \varphi(\sigma) \|_{2p}^2
+ \| \psi (\sigma)\|_{2p}^2 + \| w (\sigma)\|_{2p}^2 \big) $$
where
$$
\gamma_T = \frac{2(C+C_{B})}{T-2(C-C_{B})} < 1.
$$
Then a standard argument, similar to the one employed 
in Lemma \ref{lemma-AbsobSet}, shows that 
there exists $\gamma_{B,T},\alpha_{B,T},C_{B,T}>0$ 
such that 
$$
E_{\ell}(t) \leq \gamma_{B,T} E_{\ell}(0)e^{-\alpha_{B,T} t} + C_{B,T} \sup_{\sigma \in [0,t]}\big( \| \varphi(\sigma) \|_{2p}^2
+ \| \psi (\sigma)\|_{2p}^2 + \| w (\sigma)\|_{2p}^2 \big) .
$$
Since $T>0$ is a fixed time-step which depends on $B$, we can simply write $\gamma_{B},\alpha_{B},C_{B}$, and therefore (\ref{q-stable}) holds. \qed

\paragraph{Proof of Theorem \ref{theo-attractor2} 
(Fractal dimension).} 
We begin by observing that from the variation of parameters formula $(\ref{SolFormula})$ we obtain inequality $(\ref{local-lips})$. 
On the other hand, Lemma \ref{lemma-QEI} implies that for any bounded positively invariant set $B$, the condition 
$(\ref{StabIne})$ is valid with  
$X = H^{1}_{0}(0,L)^{3}$, $Y = L^{2}(0,L)^{3}$, $b(t) =\gamma_{B} e^{-\alpha_{B} t}$, $c(t) = C_{B}$, and 
$$
[ (\varphi,\psi,w) ]_{X} = \sqrt{ \|\varphi\|_{2p}^2+\|\psi\|_{2p}^2+\|w\|_{2p}^2} \, ,
$$
which is a compact seminorm in $X$. Then, in particular, $(\mathcal{H},S_{\ell}(t))$ is quasi-stable on the attractor $\mathbf{A}_{\ell}$. 
Therefore this attractor has finite fractal dimension from Theorem \ref{theo-A4}. \qed

\paragraph{Proof of Theorem \ref{theo-attractor2} (Regularity).} Since we know that the system is quasi-stable, 
Theorem \ref{theo-A6}, implies that any full trajectory 
$(\varphi(t),  \psi(t),  w(t),\varphi_t(t),\psi_t(t) , w_t(t))$ 
inside the attractor has regularity  
$$
\varphi_t, \psi_t, w_t \in L^{\infty}(\mathbb{R}, H^{1}_{0}(0,L))
\cap C(\mathbb{R}, L^2(0,L)) \;\; \mbox{and} \;\; \varphi_{tt}, \psi_{tt}, w_{tt} \in L^{\infty}(\mathbb{R}, L^{2}(0,L)). 
$$
Now, by continuity of the nonlinear terms, we have  
\begin{align*}
k\varphi_{xx} & = \rho_{1}\varphi_{tt} - k(\psi+\ell w)_{x} - k_{0} \ell (w_{x}-\ell\varphi)+g_{1}(\varphi_t) + f_1(\varphi,\psi,w) \in L^{\infty}(\mathbb{R};L^2(0,L)), \\
b\psi_{xx}  & = \rho_{2}\psi_{tt} +k(\varphi_{x} + \psi+\ell w) + g_{2}(\psi_t) + f_2(\varphi,\psi,w)  \in L^{\infty}(\mathbb{R};L^2(0,L)),\\
k_{0}w_{xx}  & = \rho_{1}w_{tt}+ k_{0}\ell\varphi_{x} + k\ell(\varphi_{x}+\psi+\ell w) + g_{3}(w_t)  + f_3(\varphi,\psi,w) \in L^{\infty}(\mathbb{R};L^2(0,L)).
\end{align*}
Then, elliptic regularity implies that 
$\varphi, \psi, w \in L^{\infty}(\mathbb{R}, H^2(0,L)\cap H^{1}_{0}(0,L))$, and therefore estimate (\ref{reg100}) is verified. \qed 

\paragraph{Proof of Theorem \ref{theo-attractor2} (Exponential attractors).} 
Let $\mathcal{B}$ be the bounded absorbing set of $(\mathcal{H},S_{\ell}(t))$ given by Lemma \ref{lemma-AbsobSet}.   
Then for $T>0$ and $y=(\varphi,\psi,w,\tilde{\varphi},\tilde{\psi},\tilde{w})\in \mathcal{B}$, there exists $C_{\mathcal{B}}>0$ 
such that 
$$
\|S_{\ell}(t)y \|_{\mathcal{H}}\leq C_{\mathcal{B}}, \quad 0 \leq  t \leq T . 
$$
Using this in (\ref{P01})-(\ref{P03}) we obtain   $(\varphi_{tt},\psi_{tt},w_{tt}) \in H^{-1}(0,L)^{3}$. Taking 
a larger $C_{\mathcal{B}}$ if necessary, we have  
$$
\Big\|\frac{d}{dt} S_{\ell}(t)y \Big\|_{\widetilde{\mathcal{H}}_{-1}} \leq C_{\mathcal{B}}, \quad 0 \leq  t \leq T.
$$
Consequently we obtain   
\begin{equation} \label{holder2}
\|S_{\ell}(t_1)y - S_{\ell}(t_2) y \|_{\widetilde{\mathcal{H}}_{-1}} \leq
\int_{t_1}^{t_2} \Big\|\frac{d}{dt} S_{\ell}(t)y \Big\|_{\widetilde{\mathcal{H}}_{-1}} dt \leq C_{\mathcal{B}} |t_1-t_2|, \quad  0\leq t_1 < t_2 \leq T.
\end{equation}
That is, the map $t \mapsto S_{\ell}(t)y$ is H\"older
continuous from $[0,T]$ to $\widetilde{\mathcal{H}}_{-1}$ 
(with exponent $1$). 
Therefore Theorem \ref{theo-A5} implies the existence
of a generalized exponential attractor 
$\mathbf{A}_{\ell}^{\rm exp}$ whose fractal dimension is finite in $\widetilde{\mathcal{H}}_{-1}$.

We can choose smaller extended spaces. Indeed, since 
$\widetilde{\mathcal{H}}_0 \subset \widetilde{\mathcal{H}}_{-1}$ continuously, 
given $\eta \in (0,1)$, the interpolation theorem implies that
$$
\| y \|_{\widetilde{\mathcal{H}}_{-\eta}} \leq C \,
\| y \|_{\widetilde{\mathcal{H}}_{0}}^{1-\eta} \,
\| y \|_{\widetilde{\mathcal{H}}_{-1}}^{\eta} \leq C_{\mathcal{B}}^{1-\eta} \, \| y \|_{\widetilde{\mathcal{H}}_{-1}}^{\eta}.
$$
In particular, 
$$
\| S_{\ell}(t_1)y - S_{\ell}(t_2)y \|_{\widetilde{\mathcal{H}}_{-\eta}} 
\leq C_{\mathcal{B}}^{1-\eta} \, 
\| S_{\ell}(t_1)y - S_{\ell}(t_2)y \|_{\widetilde{\mathcal{H}}_{-1}}^{\eta}.
$$
Then, combining this with (\ref{holder2}) 
we find that
$$
\|S_{\ell}(t_1)y - S_{\ell}(t_2) y \|_{\widetilde{\mathcal{H}}_{-\eta}} \leq C_{\mathcal{B}} |t_1-t_2|^{\eta}, \quad 0\leq t_1 < t_2 \leq T.
$$
This shows that $t \mapsto S_{\ell}(t)y$ is H\"older continuous in the space $\widetilde{\mathcal{H}}_{-\eta}$. 
Then the existence of a generalized exponential attractor,
with finite fractal dimension in $\widetilde{\mathcal{H}}_{-\eta}$, follows from Theorem \ref{theo-A5}. \qed 


\section{Upper-semicontinuity} 

Our last result is concerned with the convergence of attractors of the Bresse system ($\mathbf{A}_{\ell}$) to that of the Timoshenko system ($\mathbf{A}_{0}$). 
In a first understanding, we could consider the solutions of the 
Timoshenko as 
$$
(\varphi, \psi, 0 , \varphi_t ,  \psi_t , 0) \in \mathcal{H} 
\;\; \mbox{and} \;\; \mathbf{A}_{0} \subset \mathcal{H}. 
$$
Then we compare $\mathbf{A}_{\ell}$ with $\mathbf{A}_{0}$ in $\mathcal{H}$, 
as $\ell \to 0$. However, as mentioned early, in the limit the Bresse system uncouples 
into the Timoshenko system and an independent wave equation in the variable $w$.   
The model does not assert whether $w=0$. Therefore, instead extending the attractor $\mathbf{A}_{0}$ to $\mathcal{H}$, we project the attractors $\mathbf{A}_{\ell}$ onto $\mathcal{H}_{0}$. 

\begin{theorem} \label{theo-upper} Under the hypotheses of Theorem \ref{theo-attractor2}, assume further  
that $f_i$ satisfy condition $(\ref{singular50})$. Then the attractor $\mathbf{A}_{\ell}$ is upper-semicontinuous with respect 
to $\ell \to 0$, in the sense that, 
\begin{equation} \label{US1}
\lim_{\ell \rightarrow 0}{d}^{\mathcal{H}_{0}}
\bigl(\mathcal{P}\mathbf{A}_{\ell},\mathbf{A}_{0}\bigl)=0,
\end{equation}
where $d^{\mathcal{H}_{0}}$ denotes Hausdorff semi-distance, and     
$\mathcal{P}: \mathcal{H} \rightarrow \mathcal{H}_{0}$
is the projection map defined by
$ \mathcal{P}(\varphi, \psi, w,\tilde{\varphi},\tilde{\psi}, \tilde{w})=(\varphi, \psi,\tilde{\varphi}, \tilde{\psi})$. 
\end{theorem}

\proof The proof is based on the arguments in \cite{hale-raugel} and also in \cite{I2}. 
Suppose by contradiction that the statement (\ref{US1}) is false.  
Then there exists $\epsilon>0$ and a sequence $\ell_n \to 0$ 
such that 
$$
\sup_{y \in \mathbf{A}_{\ell_n}}\inf_{z\in \mathbf{A}_{0}} \Vert \mathcal{P}y - z \Vert_{\mathcal{H}_0} \ge \epsilon, 
\quad \forall \, n .
$$
Since $\mathbf{A}_{\ell_n}$ and $\mathbf{A}_{0}$ 
are compact sets, there exist $y_0^n \in \mathbf{A}_{\ell_n}$  such that 
\begin{equation} \label{US2}
\inf_{z \in \mathbf{A}_{0}} 
\Vert \mathcal{P} y_0^{n} - z \Vert_{\mathcal{H}_{0}} \geq \epsilon, \quad \forall \, n . 
\end{equation}
Let $y^{n}(t)$ the full trajectory in $\mathbf{A}_{\ell_n}$ defined by
$$
y^n(t)=(\varphi^{n}(t),\psi^n(t), w^n(t), \varphi_{t}^n(t), \psi_{t}^n(t), w_{t}^n(t)), \quad y^n(0) = y_0^n .
$$
We can assume $\ell_n \in (0,\pi/2L)$,  
and then the absorbing ball $\mathcal{B}= \overline{B}(0,{R}_0)$ given by Lemma \ref{lemma-AbsobSet} is independent of $\ell_n$. 
Then  
\begin{equation} \label{ub100}
\Vert (\varphi^n(t), \psi^n(t), w^n(t)) \Vert_{(H^{1}_0)^3}^2  + 
\Vert (\varphi_{t}^n(t), \psi_{t}^n(t), w_{t}^n(t)) \Vert_{(L^2)^3}^2  \le R_0^2, \quad \forall \, t, n. 
\end{equation}
In addition, since (\ref{norma-ok}) and (\ref{energia-domina}) are now independent of $\ell_n$, we see that 
coefficients in (\ref{q-stable}) do not depend on $\ell_n$. Then (\ref{t-regular2}) asserts that
there is $R_1>0$ 
such that  
\begin{equation} \label{ub101}
\Vert (\varphi_{t}^n(t), \psi_{t}^n(t), w_{t}^n(t))\Vert_{(H^{1}_0)^3}^2 + 
\Vert (\varphi_{tt}^n(t), \psi_{tt}^n(t), w_{tt}^n(t)) \Vert_{(L^2)^3}^2  \le R_1^2, \quad \forall \, t, n. 
\end{equation}
As in the proof Theorem \ref{theo-attractor2} (regularity), using elliptic regularity combined with $(\ref{ub100})$-$(\ref{ub101})$, 
we obtain $R_2>0$ such that
$$
\Vert (\varphi^n(t), \psi^n(t), w^n(t))\Vert_{(H^{2})^3}^2 \le R_2^2, \quad \forall \, t,n. 
$$
Consequently, 
$$
\{y^n \} \;\; \mbox{is bounded in} \;\; L^{\infty}(\mathbb{R}, H^2(0,L)^3 \times H^1(0,L)^3),
$$
$$
\{y^n_t\} \;\; \mbox{is bounded in} \;\;  L^{\infty}(\mathbb{R}, \mathcal{H}),
$$
and for every $T>0$, we have  
$$
\{y^{n}\} \;\; \mbox{is precompact in} \;\; 
C([-T,T],\mathcal{H}). 
$$
From this, there exists a subsequence $\{y^{n_k}\}$ and 
$y \in C([-T,T],\mathcal{H})$ such that   
$$
\lim_{k \rightarrow \infty} \sup_{t \in [-T,T]} \|y^{n_{k}}(t)-y(t)\|_{\mathcal{H}}=0.
$$
In particular, denoting $\tilde{z}=\mathcal{P}y$, we have 

\begin{equation} \label{US3}
\lim_{k \to \infty}\|\mathcal{P}y_0^{n_{k}} - \tilde{z}(0)\|_{\mathcal{H}_{0}} 
\to 0. 
\end{equation}
Let us show that $\tilde{z}(0)\in \mathbf{A}_0$. Indeed, 
the same argument used in Theorem \ref{theo-singular} proves that $\tilde{z}=\mathcal{P}y$ 
is a solution of the Timoshenko system (\ref{TS101})-(\ref{TS103}) in $[-T,T]$. 
Since $T>0$ is arbitrary, it follows from (\ref{ub100}) that 
$\tilde{z}(t)$ is a bounded full trajectory for the Timoshenko system and thus  
$\tilde{z}(0) \in \mathbf{A}_0$. Therefore, (\ref{US3}) contradicts (\ref{US2}). 
This completes the proof of the Theorem \ref{theo-upper}. \qed



\appendix
\section*{Appendix}
\setcounter{equation}{0}
\setcounter{theorem}{0}
\renewcommand{\thetheorem}{A.\arabic{theorem}}
\renewcommand{\theequation}{A.\arabic{equation}}
\renewcommand{\thesection}{A}

In this appendix we have collected some useful results concerning attractors of nonlinear 
infinite dimensional dynamical systems. Most of results can be found in classical references such as
\cite{babin-visik,hale,lady,temam}. We shall follow more closely \cite[Chapter 7]{Yellow}.

\medskip
\noindent{\bf Definitions.}
\medskip

A dynamical system is a pair $(H,S(t))$ where $H$ is a complete metric space and $S(t)$ 
is a strongly continuous semigroup of $H$. 
Then, a global attractor for $(H,S(t))$ is a compact set $\mathbf{A} \subset H$ that is fully
invariant and uniformly attracting, that is,
$$
S(t)\mathbf{A} = \mathbf{A} \;\; \mbox{and} \;\; \lim_{t\to \infty} d^{H}(S(t)B,\mathbf{A})=0,
$$
for any bounded set $B\subset H$, where $d^H$ denotes the Hausdorff semidistance 
$$
d^H(A,B) = \sup_{a\in A} \inf_{b \in B} d(a,b). 
$$
Roughly speaking, the existence of a global attractor relies on two properties, namely, dissipativeness and compactness.
A dynamical
system is called dissipative if it has an absorbing set, that is, a bounded set $\mathcal{B}\subset H$ that attracts any bounded set $B$ in a finite time $T_B>0$. In other words,
$$
S(t)B \subset \mathcal{B}, \quad t \ge T_B.
$$
As for compactness, a dynamical system $(H,S(t))$ is called
asymptotically compact if for any bounded set $B\subset H$, and sequence
$\{x_{k}\} \subset B$, the sequence $\{S(t_k)x_{k}\}$ has convergent subsequence
whenever $t_k \to \infty$. This condition is often hard to prove and
some compactness criteria are used instead. The fractal dimension of a compact set $A \subset H$ is a number defined by
$$
\dim_F A = \limsup_{\varepsilon \to 0} \frac{\ln N_{\varepsilon}(A)}{\ln(1/{\varepsilon})} ,
$$
where $N_{\varepsilon}(A)$ is the minimal quantity of closed balls of radius $2 \varepsilon$ 
necessary to cover $A$. A compact set $\mathbf{A}^{\rm exp} \subset H$ is called a generalized exponential attractor if: 

$(i)$ it is positively invariant, 

$(ii)$ it attracts exponentially
fast the trajectories from any bounded set of $H$, 

$(iii)$ it has finite fractal dimension 
in an extended space $\widetilde{H} \supseteq {H}$. 

\smallskip

In what follows we present some abstract results that 
guarantee the existence of smooth attractors with finite fractal dimension.  

\medskip
\noindent{\bf A compactness criterion.}
\medskip

The following criterion for asymptotic compactness is very useful
for hyperbolic like systems. It was presented in
\cite{LasieckaI,Yellow} and involves a function
$\Psi: B\times B \to \mathbb{R}$ such that
\begin{equation} \label{contractive-f}
\liminf_{m\rightarrow \infty}\liminf_{n\rightarrow \infty}\Psi(y_{n},y_{m})=0 ,
\end{equation}
for every sequence $\{y_{n}\} \subset B$, where $B$ is a bounded set of $H$.

\begin{theorem}
{\rm \cite[Theorem 7.1.11]{Yellow}}
\label{theo-A1} Let $({H}, S(t))$
be a dynamical system where $H$ is a Banach space.
Assume that for any bounded positively invariant set $B \subset H$ and any $\epsilon> 0$,
there exists a time $T =T_{\epsilon,B}$ and a function
$\Psi_{\epsilon, B, T}: B\times B \to \mathbb{R}$ satisfying $(\ref{contractive-f})$ such that
$$
\|S(T)y^{1}-S(T)y^{2} \|_{{H}} \leq \epsilon
+ \Psi_{\epsilon, B, T}(y^{1},y^{2}), \quad \forall \: y^{1},y^{2} \in B.
$$
Then $(H, S(t))$ is asymptotically compact.
\end{theorem}


\medskip
\noindent{\bf Gradient systems.}
\medskip

Gradient systems have more specialized
dissipativeness because they admit a strict Lyapunov function. More precisely, a functional $\Phi:H \to \mathbb{R}$ is
a strict Lyapunov function for a system $({H},S(t))$ if,

$(i)$ the map $t \mapsto \Phi(S(t)z)$ is non-increasing for any $z \in H$,

$(ii)$ if $\Phi(S(t)z)=\Phi(z)$ for all $t$, then $z$ is a stationary point of $S(t)$.

\noindent Attractors of gradient systems may have further geometric properties. Let $\mathcal{N}$ 
be the set of stationary points of $S(t)$. Then the unstable 
manifold $\mathbb{M}_{+}(\mathcal{N})$ is the family of $y \in H$  
such that there exists a full trajectory $u(t)$ satisfying  
$$
u(0) = y \;\; \mbox{and} \;\; \lim_{t\to -\infty} dist (u(t), \mathcal{N}) = 0.
$$
The following result is well-known. See for instance \cite[Corollary 7.5.7]{Yellow}.

\begin{theorem} \label{theo-A2}
Let $({H},S(t))$ be an asymptotically compact gradient system
with the corresponding Lyapunov functional denoted by $\Phi$. Suppose that
\begin{equation} \label{SSE}
\Phi(z) \to \infty \;\; \mbox{if and only if} \;\; \Vert z \Vert_{H} \to \infty,
\end{equation}
and that the set of stationary points $\mathcal{N}$ is bounded.
Then $({H},S(t))$ has a compact global attractor which
coincides with the unstable manifold $\mathbb{M}_{+}(\mathcal{N})$.
\end{theorem}


\medskip
\noindent{\bf Quasi-stable systems.}
\medskip

Let $X,Y$ be two reflexive Banach spaces with $X$ compactly embedded into $Y$ and
put ${H}= X \times Y$. Consider the dynamical system $({H},S(t))$ given by
\begin{equation} \label{DSystem}
S(t)y= (u(t),u_{t}(t)), \qquad y=(u(0), u_{t}(0)) \in {H},
\end{equation}
where 
the functions $u$ have regularity
\begin{equation}\label{Regularity}
u \in C([0,\infty); X)\cap C^{1}([0,\infty); Y).
\end{equation}
Then we say  
that it is quasi-stable on a set
$B \subset {H}$, if there exist a compact semi-norm $[ \, \cdot \, ]_{X}$
on $X$ and nonnegative scalar functions
$a(t)$ and $c(t)$, locally bounded in $[0,\infty)$, and
$b(t)\in L^1(0,\infty)$ with $\lim_{t\to \infty}b(t)=0$, such that,
\begin{equation} \label{local-lips}
\| S(t)y^1 - S(t)y^2 \|_{{H}}^2 \leq a(t)\|y^1 - y^2 \|_{{H}}^2 ,
\end{equation}
and
\begin{equation} \label{StabIne}
\|S(t)y^1 - S(t)y^2 \|_{{H}}^2 \leq  b(t) \| y^1 -y^2 \|_{{H}}^2 
+ c(t) \sup_{0<s<t} [ u^1(s)-u^2(s) ]_{X}^2 ,
\end{equation}
for any $y^1,y^2 \in B$, where $S(t)y^{i}= (u^{i}(t),u_{t}^{i}(t))$, $i=1,2$.

\smallskip

The first property of quasi-stable system is the asymptotic compactness.

\begin{theorem}
{\rm \cite[Proposition 7.9.4]{Yellow}}
\label{theo-A3}
Let $({H},S(t))$ be a dynamical system
given by $(\ref{DSystem})$ and satisfying $(\ref{Regularity})$.
Suppose that the system is quasi-stable on every bounded positively
invariant set $B$ of ${H}$.
Then $({H},S(t))$ is asymptotically compact.
\end{theorem}


\medskip
\noindent{\bf Fractal dimension and exponential attractors.}
\medskip

Quasistability also implies that global attractors have finite fractal-dimension. 

\begin{theorem}
{\rm \cite[Theorem 7.9.6]{Yellow}}
\label{theo-A4}
Let $({H},S(t))$ be a dynamical system
given by $(\ref{DSystem})$ and satisfying $(\ref{Regularity})$.
Suppose that it has a global attractor $\mathbf{A}$ and it is
quasi-stable on it. Then $\mathbf{A}$ has finite fractal dimension.
\end{theorem}

\begin{theorem}
{\rm \cite[Theorem 7.9.9]{Yellow}}
\label{theo-A5}
Let $({H},S(t))$ be a dissipative dynamical system
satisfying $(\ref{DSystem})$-$(\ref{Regularity})$ and quasi-stable on some bounded absorbing set $\mathcal{B}$.
In addition, suppose there exists an extended space $\widetilde{H} \supseteq {H}$ such that,
for each $T>0$,
$$ 
\Vert S(t_1)y - S(t_2)y \Vert_{\widetilde{H}} \leq C_{\mathcal{B}T} |t_1 - t_2|^{\eta},
\qquad t_1,t_2 \in [0,T], \;\; y \in \mathcal{B},
$$ 
where $C_{\mathcal{B}T}>0$ and $\eta \in (0,1]$ are constants.
Then this system has a generalized exponential attractor $\mathbf{A}^{\rm exp} \subset H$
with finite fractal dimension in $\widetilde{H}$.
\end{theorem}


\medskip
\noindent{\bf Regularity of attractors.}
\medskip

In many cases a quasi-stable system
has global attractor more regular than its phase space. The next result
is about gain of regularity in the $t$ variable. 

\begin{theorem}
{\rm \cite[Theorem 7.9.8]{Yellow}}
\label{theo-A6}
Let $({H},S(t))$ be a dynamical system satisfying $(\ref{DSystem})$-$(\ref{StabIne})$ with 
$c(t)$ bounded. Suppose that it has a global attractor $\mathbf{A}$ and it is
quasi-stable on it. Then any full trajectory $(u(t),u_t(t))$ in the attractor 
have additional regularity 
$$
u_t \in L^{\infty} (\mathbb{R}, X) \cap 
C(\mathbb{R}, Y) \;\; \mbox{and} \;\; 
u_{tt} \in L^{\infty} (\mathbb{R}, Y).
$$
In addition, 
\begin{equation} \label{t-regular2}
\Vert u_t(t)\Vert_X^2 + \Vert u_{tt}(t)\Vert_Y^2 \le R^2, \quad t \in \mathbb{R},
\end{equation}
where $R>0$ depends on $\sup_{t>0} c(t)$, $\mu_X$,  
and on the embedding $X \hookrightarrow Y$. 
\end{theorem}
The corresponding gain of regularity in the $x$ variable is usually obtained from elliptic regularity.

\bigskip

\paragraph{Acknowledgment.} The first author was partially supported by CNPq Grant 310041/2015-5. 
The second author was supported by CAPES/PDSE Grant 04724/2014.

\end{document}